\documentclass[aap,MSNbibl,seceqn,citesort,dvips]{arximspdf}
\usepackage{mathbh}

\doi{10.1214/11-AAP792}
\volume{22}
\issue{3}
\pubyear{2012}
\firstpage{1167}
\lastpage{1214}

\makeatletter

\newcommand{\bb}{\mathbb}

\newcommand{\PP}{{\bb P}}
\newcommand{\E}{{\bb E}}
\newcommand{\real}{{\bb{R}}}
\newcommand{\definedas}{\stackrel{\Delta}{=}}
\newcommand{\equalinlaw}{\stackrel{\mathcal L}{=}}

\newcommand{\Cov}{\operatorname{Cov}}
\newcommand{\indic}{{\mathbh1}}

\newtheorem{theorem}{Theorem}[section]
\newtheorem{lemma}[theorem]{Lemma}
\newtheorem{corollary}[theorem]{Corollary}
\newtheorem{proposition}[theorem]{Proposition}

\newproclaim{definition}[theorem]{Definition}
\newproclaim{remark}[theorem]{Remark}
\newproclaim{Example}[theorem]{Example}
\newproclaim{algorithm}[theorem]{Algorithm}

\setattribute{abstract}{skip}{25\p@}

\makeatother

\begin{document}
\begin{frontmatter}

\title{Efficient Monte Carlo for high excursions of~Gaussian random fields}
\runtitle{Monte Carlo for random fields}

\begin{aug}
\author[A]{\fnms{Robert J.} \snm{Adler}\thanksref{t1}\ead[label=e1]{robert@ee.technion.ac.il}\ead[label=u1,url]{webee.technion.ac.il/people/adler}},
\author[B]{\fnms{Jose H.} \snm{Blanchet}\thanksref{t2}\ead[label=e2]{jose.blanchet@columbia.edu}\ead[label=u2,url]{www.ieor.columbia.edu/fac-bios/blanchet/faculty.html}} and
\author[C]{\fnms{Jingchen} \snm{Liu}\corref{}\thanksref{t3}\ead[label=e3]{jcliu@stat.columbia.edu}\ead[label=u3,url]{stat.columbia.edu/\textasciitilde jcliu}}

\runauthor{R. J. Adler, J. H. Blanchet and J. Liu}
\affiliation{Technion, Columbia Uninversity and Columbia Uninversity}
\address[A]{R. J. Adler\\
Faculty of Electrical Engineering\\
Technion\\
Haifa\\
Israel 32000 \\
\printead{e1}\\
\printead{u1}}
\address[B]{J. H. Blanchet\\
Department of Industrial Engineering\\
\quad and Operations Research\\
Columbia Uninversity\\
340 S.W. Mudd Building\\
500 W. 120 Street\\
New York, New York 10027\\
USA\\
\printead{e2}\\
\printead{u2}}
\address[C]{J. Liu\\
Department of Statistics\\
Columbia University\\
1255 Amsterdam Avenue Room 1030\\
New York, New York 10027\\
USA\\
\printead{e3}\\
\printead{u3}}
\end{aug}

\thankstext{t1}{Supported in part by US--Israel Binational
Science Foundation, 2008262, ISF 853/10 and NSF Grant DMS-08-52227.}

\thankstext{t2}{Supported in part by NSF Grant DMS-09-02075, CMMI-0846816
and CMMI-1069064.}

\thankstext{t3}{Supported in part by Institute of Education Sciences,
U.S. Department of Education, through Grant R305D100017 and
NSF CMMI-1069064.}

% HISTORY:
\received{\smonth{5} \syear{2010}}
\revised{\smonth{3} \syear{2011}}

% ABSTRACT
%
\begin{abstract}
Our focus is on the design and analysis of efficient Monte Carlo
methods for computing tail probabilities for the suprema of Gaussian
random fields, along with conditional expectations of functionals of
the fields given the existence of excursions above high levels, $b$.
Na\"{\i}ve Monte Carlo takes an exponential, in $b$, computational
cost to estimate these probabilities and conditional expectations for
a prescribed relative accuracy. In contrast, our Monte Carlo
procedures achieve, at worst, polynomial complexity in $b$, assuming
only that the mean and covariance functions are H\"{o}lder continuous.
We also explain how to fine tune the construction of our procedures in
the presence of additional regularity, such as homogeneity and
smoothness, in order to further improve the efficiency.
\end{abstract}

% KEYWORDS
%

\end{frontmatter}

%s1 ###
\section{Introduction}

This paper centers on the design and analysis of
efficient Monte Carlo techniques for computing probabilities and conditional
expectations related to high excursions of Gaussian random fields. More
specifically, suppose that $f\dvtx T\times\Omega\rightarrow\real$
is a continuous Gaussian random field over a $d$-dimensional
compact set $T\subset\real^{d}$. (Additional
regularity conditions on $T$
will be imposed below, as needed.)

Our focus is on tail probabilities of the form
%
%e1.1 ###
%
\begin{equation}
\label{tailprobequn1} \label{TP}%
w(b)=\PP\Bigl(\max_{t\in T}f( t) >b\Bigr)
\end{equation}
and on conditional expectations %
%
%e1.2 ###
%
\begin{equation}\label{CE}%
\E\Bigl[\Gamma(f)\big|\max_{t\in T}f( t) >b\Bigr]
\end{equation}
as $b\rightarrow\infty$, where $\Gamma$ is a functional
of the field, which, for concreteness we take to be positive and bounded.

While much of the paper will concentrate on estimating the exceedance
probability~(\ref{tailprobequn1}), it is important to note that our
methods, based on importance sampling, are broadly
applicable to the efficient evaluation of conditional expectations of the
form~(\ref{CE}). Indeed, as we shall explain
at the end of Section~\ref{SecIS}, our approach to
efficient importance sampling is based on a procedure which
mimics the conditional distribution of $f$, given that
$\max_{T}f(t) >b$. Moreover, once an efficient (in a
precise mathematical sense described in Section~\ref{SecSet}) importance
sampling procedure is in place, it follows under mild regularity
conditions on
$\Gamma$ that an efficient estimator for~(\ref{CE}) is
immediately obtained by exploiting an efficient estimator for~(\ref{TP}).

The need for an efficient estimator of $\omega(b)$ should be reasonably
clear. Suppose one could simulate
\[
f^{\ast}\triangleq\sup_{t\in T}f( t)
\]
exactly\vspace*{1pt} (i.e., without bias) via na\"ive Monte Carlo. Such an approach
would typically require a number of replications of $f^{\ast}$ which
would be
exponential in $b^{2}$ to obtain an accurate estimate (in relative terms).
Indeed, since in great generality (see~\cite{SheppLan70})
$w( b) =\exp( -cb^{2}+o( b^{2}))$ as $b\rightarrow
\infty$ for some $c\in(0,\infty)$, it follows that the
average of $n$ i.i.d. Bernoulli trials each with success parameter\vspace*{1pt} $w(
b) $
estimates $w( b) $ with a relative mean squared error equal to $n^{-1/2}(1-w(
b) )^{1/2}/w( b) ^{1/2}$. To control the size of the error therefore
requires\setcounter{footnote}{3}\footnote{Given $h$ and $g$ positive, we shall use the
familiar asymptotic
notation $h(x)=O( g( x) ) $ if there is $c<\infty$ such that $h( x)
\leq cg(
x) $ for all $x$ large enough; $h( x) =\Omega( g( x) ) $ if $h( x)
\geq cg( x)
$ if $x$ is sufficiently large and $h( x) =o( g( x) ) $ as
$x\rightarrow
\infty$ if $h( x) /g( x) \rightarrow0$ as $x\rightarrow\infty$; and
$h(x) =
\Theta(g(x))$ if $h(x)= O(g(x))$ and $h(x)= \Omega(g(x))$.}
$n=\Omega(w( b) ^{-1})$, which is typically prohibitively large. In addition,
there is also a problem in that typically $f^{\ast}$ cannot be
simulated exactly and that some discretization of $f $ is required.

Our goal is to introduce and analyze simulation estimators that can be
applied to a general class of Gaussian fields and that can be shown to require
at most a polynomial number of function evaluations in $b$ to obtain estimates
with a prescribed relative error. The model of computation that we use to
count function evaluations and the precise definition of an algorithm with
polynomial complexity is given in Section~\ref{SecSet}. Our proposed
estimators are, in particular, \textit{asymptotically optimal}. (This
property, which is a popular notion in the context of
rare-event simulation (cf. \mbox{\cite{ASGLYNN,BUCK04}}) essentially requires
that the second moments of estimators decay\vadjust{\goodbreak}
at the same exponential rate as the square of the
first moments.) The polynomial complexity of our estimators requires to
assume no more than that the underlying
Gaussian field is H\"{o}lder continuous; see Theorem~\ref{ThmGeneral} in
Section~\ref{SecMain}. Therefore, our methods provide means for efficiently
computing probabilities and expectations
associated with high excursions of Gaussian random fields in wide
generality.\looseness=1

In the presence of enough smoothness, we shall also show
how to design estimators that
can be shown to be \textit{strongly efficient}, in the sense that their
associated coefficient of variation remains uniformly bounded as
\mbox{$b\rightarrow\infty$}. Moreover, the associated path
generation of the conditional field (given a~high excursion) can,
basically, be carried out
with the same computational complexity as the unconditional sampling
procedure (uniformly in $b$). This is Theorem~\ref{ThmHom} in Section
\ref{SecMain}.

High excursions of Gaussian random fields appear in wide number of applications,
including, but not limited to:
\begin{itemize}
\item Physical oceanography: Here the random field can be water
pressure or surface temperature. See~\cite{ADLER-OCEAN} for many examples.
\item Cosmology: This includes the analysis of COBE and WMAP microwave data
on a sphere or galactic density data; for example, \cite
{Bardeen,Shandarin1,Shandarin2}.
\item Quantum chaos: Here random planar waves replace deterministic
(but unobtainable) solutions of Schrodinger equations; for example, the
recent review
\cite{DENNIS-REV}.
\item Brain mapping: This application is
the most developed and very widely used. For example,
the paper by Friston et al.~\cite{FRISTON94}
that introduced random field methodology to the brain imaging
community has, at the time of writing, over 4,500 (Google) citations.
\end{itemize}

Many of these applications deal with twice differentiable,
constant variance random
fields, or random fields that have been normalized
to have constant variance, the reason being that they require estimates
of the tail probabilities~(\ref{tailprobequn1}) and these are only really
well known in the smooth,
constant (unit) variance case. In particular, it is known that, with
enough smoothness
assumptions,
%
%e1.3 ###
%
\begin{eqnarray}
\label{mainresultintro}
&&\liminf_{b \rightarrow\infty} -b^{-2} \log\Bigl| P \Bigl( \sup_{t \in
T} f(t) \geq b \Bigr) - E\bigl(\chi\bigl(\{t\in T\dvtx f(t)\geq b\} \bigr) \bigr)
\Bigr|\nonumber\\[-8pt]\\[-8pt]
&&\qquad\geq\frac{1}{2} + \frac{1}{2\sigma^2_c},\nonumber
\end{eqnarray}
where $\chi(A)$ is the Euler characteristic of the set $A$, and
the term $\sigma^2_c$ is related to a geometric quantity known as the
critical radius of $T$ and depends on the covariance structure of $f$;
see~\cite{AdlerTaylor07,TaylorTakemuraAdler2003Validity}
for details.
Since both the probability and the expectation in~(\ref{mainresultintro})
are typically $O(b^\ell\exp(-b^2/2))$ for some \mbox{$\ell\geq0$} and large\vadjust{\goodbreak}
$b$, a more user friendly (albeit not quite as correct) way to
write~(\ref{mainresultintro}) is
%
%e1.4 ###
%
\begin{equation}
\label{mainresultintroeasier}
P \Bigl( \sup_{t \in T} f(t) \geq b \Bigr) \approx
E\bigl(\chi\bigl(\{t\in T\dvtx f(t)\geq b\} \bigr) \bigr)
\times\bigl(1+ O(e^{-cb^2})\bigr)
\end{equation}
for some $c$.
%%\liminf_{b \rightarrow\infty} -b^{-2} \log| P ( \sup_{t \in T}
%f(t) \geq b ) - E(\chi(\{t\in T\dvtx f(t)\geq b\} ) ) | \geq
%| P ( \sup_{t \in T} f(t) \geq b ) - E(\chi(\{t\in T\dvtx f(t)
%P(\sup_{t\in T}f(t) > b )\times e^{-\frac{b^2}{2\sigma_c^2}
%+o(b^2)},

The expectation in~(\ref{mainresultintro}) and
(\ref{mainresultintroeasier})
has an explicit form that is
readily computed for Gaussian and Gaussian-related random fields of
constant variance (see~\cite{AdlerTaylor07,ADT10} for details),
although if
$T$ is geometrically complicated or the covariance of $f$ highly nonstationary
there
can be terms in the expectation that can only be evaluated numerically or
estimated statistically; for example,~\cite{Adler-Bartz-Kou,Taylor-Worsley-JASA}.
Nevertheless, when available,~(\ref{mainresultintro})
provides excellent approximations and simulation studies have shown
that the
approximations are numerically useful for quite moderate $b$, of the order
of 2 standard deviations.

However, as we have already noted,~(\ref{mainresultintro}) holds only
for constant variance fields, which also need to be twice differentiable.
In the case of less smooth $f$, other classes of results occur, in
which the
expansions are less reliable and, in addition, typically involve the unknown
Pickands' constants; cf.~\cite{AzWschebor09,PiterbargGaussian}.

These are some of the reasons why, despite a well developed theory,
Monte Carlo techniques still have a significant role to play in
understanding the behavior of
Gaussian random fields at high levels.
The estimators proposed in this paper basically reduce the rare event
calculations associated to high excursions in Gaussian random fields to
calculations that are roughly comparable to the evaluation of
expectations or
integrals in which no rare event is involved. In other words, the
computational complexity required to implement the estimators discussed here
is similar in some sense to the complexity required to evaluate a given
integral in finite dimension or an expectation where no tail parameter is
involved. To the best of our knowledge these types of reductions have
not been
studied in the development of numerical methods for high excursions of random
fields. This feature distinguishes the present work from the
application of
other numerical techniques that are generic (such as quasi-Monte Carlo and
other numerical integration routines) and that in particular might be also
applicable to the setting of Gaussian fields.

Contrary to our methods, which are designed to have provably good performance
uniformly over the level of excursion, a generic numerical
approximation procedure, such as quasi-Monte Carlo, will typically
require an
exponential increase in the number of function evaluations in order to
maintain a prescribed level of relative accuracy. This phenomenon is
unrelated to the setting of Gaussian random fields.
In particular, it can be easily seen to happen even when evaluating a
one-dimensional integral with a small integrand. On the other hand, we believe
that our estimators can, in practice, be easily combined with
quasi-Monte Carlo or other numerical integration methods. The rigorous
analysis of such hybrid approaches, although of great interest,
requires an extensive treatment and will be pursued in the future. As an
aside, we note that quasi-Monte Carlo techniques have
been used in the excursion analysis
of Gaussian random fields in~\cite{AzWschebor09}.

The remainder of the paper is organized as follows. In Section \ref
{SecSet} we
introduce the basic notions of polynomial algorithmic complexity, which are
borrowed from the general theory of computation. Section~\ref{SecMain}
discusses the main results in light of the complexity considerations
of Section~\ref{SecSet}. Section~\ref{SecIS} provides a brief
introduction to importance sampling, a simulation technique that we shall
use heavily in the design of our algorithms. The analysis of finite fields,
which is given in Section~\ref{SecFinite}, is helpful to develop the basic
intuition behind our procedures for the general case. Section \ref
{SecHolder} provides the
construction and analysis of a polynomial time algorithm for high excursion
probabilities of H\"{o}lder continuous fields. Finally, in
Section~\ref{SecSmooth}, we add additional smoothness assumptions along with
stationarity and explain how to fine tune the construction of our
procedures in order to further improve efficiency in these cases.

%s2 ###
\section{Basic notions of computational complexity}\label{SecSet}

In this section we shall discuss some general notions of efficiency and
computational complexity related to the approximation of the
probability of
the rare events $\{B_{b}\dvtx b\geq b_{0}\}$, for which
$P( B_{b}) \rightarrow0$ as $b\rightarrow\infty$. In
essence, efficiency means that computational complexity is, in some
sense, controllable,
uniformly in $b$. A notion that is popular in the
efficiency analysis of Monte Carlo methods for rare events is weak efficiency
(also known as asymptotic optimality) which requires that the
coefficient of
variation of a given estimator, $L_{b}$ of $P( B_{b})$,
to be dominated by $1/P(
B_{b})^{\varepsilon}$ for any $\varepsilon>0$. More formally,
we have:
\begin{definition}
\label{DefPoly} A family of estimators $\{L_{b}\dvtx b\geq b_{0}\}$
is said to be
{polynomially efficient} for estimating $P( B_{b}) $ if
$E(L_{b})=P( B_{b}) $, and there exists a~$q\in(0,\infty
)$ for
which
%
%e2.1 ###
%
\begin{equation} \label{bnd}%
\sup_{b\geq b_{0}}\frac{\operatorname{Var}(L_{b})}{[P( B_{b}
)]^{2}\vert{\log
P}( B_{b}) \vert^{q}}<\infty.
\end{equation}
Moreover, if~(\ref{bnd}) holds with $q=0$, then the family is said to be
{strongly efficient}.
\end{definition}

Below we often refer to $L_{b}$ as a \textit{strongly} (\textit{polynomially})
\textit{efficient estimator}, by
which we mean
that the family $\{L_{b}\dvtx b>0\}$ is strongly\vspace*{1pt} (polynomially) efficient.
In order to understand the nature of this definition
let $\{L_{b}^{( j) }, 1\leq j\leq n\}$ be a collection of
i.i.d.
copies of $L_{b}$. The averaged estimator
\[
\widehat{L}_{n}( b) =\frac{1}{n}\sum
_{j=1}^{n}L_{b}^{(
j) }%
\]
has a \textit{relative} mean squared error equal to
$[\operatorname{Var}( L_{b})]^{1/2}/[n^{1/2}P( B_{b}) ]$. A simple consequence of
Chebyshev's inequality is that%
\[
P\bigl( \vert\widehat{L}_{n}( b) /P(
B_{b})
-1\vert\geq\varepsilon\bigr) \leq\frac{\operatorname{Var}(
L_{b})
}{\varepsilon^{2} n P[( B_{b})] ^{2}}.
\]
Thus, if $L_{b}$ is polynomially efficient, and we wish to compute
$P(
B_{b}) $ with at most $\varepsilon$ relative error and at least
$1-\delta$ confidence, it suffices to simulate\looseness=1
\[
n=\Theta( \varepsilon
^{-2}\delta^{-1}\vert{\log P}( B_{b}) \vert
^{q})
\]\looseness=0
i.i.d. replications of $L_{b}$. In fact, in the presence of
polynomial efficiency, the bound $n=\Theta( \varepsilon
^{-2}\delta
^{-1}\vert{\log P}( B_{b}) \vert^{q})
$ can be
boosted to $n=\Theta( \varepsilon^{-2}\log(\delta^{-1})
\vert\break{\log
P}( B_{b}) \vert^{q}) $ using the so-called median
trick; see~\cite{NiePoka09}.

Naturally, the cost per replication must also be considered in the analysis,
and we shall do so, but the idea is that evaluating $P(
B_{b}) $
via crude Monte Carlo would require, given $\varepsilon$ and $\delta$,
$n=\Theta( 1/P( B_{b}) ) $ replications.
Thus a~polynomially efficiently estimator makes the evaluation of $P(
B_{b}) $ \textit{exponentially faster} relative to crude Monte Carlo,
at least in terms of the number of replications.

Note that a direct application of deterministic algorithms (such as
quasi-Monte Carlo or quadrature integration rules) might improve (under appropriate
smoothness assumptions) the computational complexity relative to Monte Carlo,
although
only by a polynomial rate (i.e., the \textit{absolute} error~decreases to zero
at rate $n^{-p}$ for $p>1/2$, where $n$ is the number of function evaluations
and $p$ depends on the dimension of the function that one is
integrating; see,
e.g.,~\cite{ASGLYNN}). We believe that the procedures that we
develop in this paper can guide the construction of efficient deterministic
algorithms with small relative error and with complexity that scales at a
polynomial rate in $\vert{\log P}( B_{b})
\vert$. This
is an interesting research topic that we plan to explore in the
future.

An issue that we shall face in designing our Monte Carlo procedure is that,
due to the fact that $f$ will have to be discretized,
the corresponding estimator~$\widetilde L_{b}$ will not be unbiased. In
turn, in order to control the
relative bias with an effort that is comparable to the bound on the
number of
replications discussed in the preceding paragraph, one must verify that the
relative bias can be reduced to an amount less than $\varepsilon$ with
probability at least $1-\delta$ at a computational cost of the form
$O(\varepsilon^{-q_{0}}\vert{\log P}( B_{b})\vert
^{q_{1}})$. If $\widetilde{L}_{b}(\varepsilon)$ can be\vspace*{1pt} generated with
$O(\varepsilon^{-q_{0}}\vert{\log P}( B_{b})\vert
^{q_{1}})$ cost, and satisfying $\vert P( B_{b})
-E\widetilde{L}_{b}(\varepsilon)\vert\leq\varepsilon P( B_{b}) $, and
if
\[
\sup_{b>0}\frac{\operatorname{Var}(\widetilde{L}_{b}(\varepsilon))}{P(
B_{b})
^{2}\vert{\log P}( B_{b}) \vert^{q}}<\infty
\]
for some $q\in(0,\infty)$, then $\widehat{L}_{n}^{\prime}(
b,\varepsilon) =\sum_{j=1}^{n}\widetilde
{L}_{b}^{(j)}(\varepsilon)/n$,
where the $\widetilde{L}_{b}^{(j)}(\varepsilon)$'s are i.i.d.
copies of
$\widetilde{L}_{b}(\varepsilon)$, satisfies
\[
P\bigl(|\widehat{L}_{n}^{\prime}( b,\varepsilon) /P(
B_{b}) -1|\geq2\varepsilon\bigr)\leq\frac{\operatorname{Var}(\widetilde{L}_{b}%
(\varepsilon))}{\varepsilon^{2}\times n\times P( B_{b}) ^{2}}.\vadjust{\goodbreak}
\]
Consequently, taking $n=\Theta( \varepsilon^{-2}\delta
^{-1}\vert{\log
P}( B_{b}) \vert^{q}) $ suffices to give an estimator
with at most $\varepsilon$ relative error and $1-\delta$ confidence,
and the
total computational cost is $\Theta( \varepsilon
^{-2-q_{0}}\delta^{-1}\vert{\log P}( B_{b})
\vert^{q+q_{1}%
}) $.

We shall measure computational cost in terms of function evaluations such
as a single addition, a multiplication, a comparison, the generation of a
single uniform random variable on $T$, the generation of a single standard
Gaussian random variable and the evaluation of $\Phi( x) $
for fixed $x\geq0$, where
\[
\Phi( x) = 1-\Psi(x) =\frac{1}{ \sqrt{2\pi}}\int_{-\infty
}^{x}e^{-s^{2}/2}
\,ds.
\]
All of these function evaluations are assumed to cost at most a fixed
amount~$\mathbf{c}$ of computer time. Moreover, we shall also assume
that first- and \mbox{second-order} moment characteristics of the
field, such as $\mu( t) =Ef(t) $ and $C( s,t) =\Cov( f( t) ,
f( s) ) $
can be computed in at most $\mathbf{c}$ units of computer time for each
$s,t\in T$. We note that similar models of computation are often used
in the complexity theory of continuous problems; see~\cite{TraubWW88}.

The previous discussion motivates the next definition which has its
roots in
the general theory of computation in both continuous and discrete settings
\cite{MitzUpf05,TraubWW88}. In particular, completely analogous notions
in the setting of complexity theory of continuous problems lead to the notion
of ``tractability'' of a computational problem~\cite{Wos96}.
\begin{definition}
A Monte Carlo procedure is said to be a \textit{fully
polynomial randomized approximation scheme} (FPRAS) for estimating
$P(
B_{b}) $ if, for some $q,q_{1},q_{2}\in\lbrack0,\infty)$, it
outputs an averaged estimator that is guaranteed to have
at most $\varepsilon>0$ relative error with confidence at least
$1-\delta
\in(0,1)$ in $\Theta(\varepsilon^{-q_{1}}\delta^{-q_{2}}|{\log
P}(
B_{b}) |^{q})$ function evaluations.
\end{definition}

The terminology adopted, namely FPRAS, is borrowed from the complexity theory
of randomized algorithms for counting~\cite{MitzUpf05}. Many counting
problems can be expressed as rare event estimation problems. In such
cases it
typically occurs that the previous definition (expressed in terms of a
rare event probability) coincides precisely with the standard counting
definition of a FPRAS (in which there is no reference to any rare event to
estimate). This connection is noted, for instance, in \cite
{Blanchet07}. Our
terminology is motivated precisely by this connection.

By letting $B_{b}=\{f^{\ast}>b\}$, the goal in this paper is to design
a class
of fully polynomial randomized approximation schemes that are
applicable to a
general class of Gaussian random fields. In turn, since our\vadjust{\goodbreak} Monte Carlo
estimators will be based on importance sampling, it turns out that we
shall also be
able to straightforwardly construct FPRASs to estimate quantities such as
$E[\Gamma(f)|{\sup_{t\in T}f}(t)>b]$ for a suitable class of functionals
$\Gamma$ for which $\Gamma( f) $ can be
computed with an error of at most $\varepsilon$ with a cost that is polynomial
as function of $\varepsilon^{-1}$. We shall discuss this observation
in Section
\ref{SecIS}, which deals with properties of importance sampling.

%s3 ###
\section{Main results}\label{SecMain}

In order to state and discuss our main results we need some notation.
For each $s,t\in T$ define%
\begin{eqnarray*}
\mu(t)&=&E(f(t)),\qquad \mathcal{C}( s,t) =\Cov(f(
s)
,f( t) ),\\
\sigma^{2}(t)&=&\mathcal{C}(t,t)>0,\qquad
r(s,t)=\frac{\mathcal{C}(s,t)}{\sigma(s)\sigma(t)}.
\end{eqnarray*}
Moreover, given $x\in\real^{d}$ and $\beta>0$ we write $\vert
x\vert
={\sum_{j=1}^{d}}\vert x_{j}\vert$, where $x_{j}$ is the $j$th
component of $x$. We shall assume that, for each fixed $s,t\in T$, both
$\mu( t) $ and $\mathcal{C}( s,t) $ can be evaluated
in at most $\mathbf{c}$ units of computing time.

Our first result shows that under modest continuity conditions on $\mu$,
$\sigma$ and~$r $ it is possible to construct an explicit
FPRAS for $w( b) $ under the following regularity conditions:
\begin{longlist}[(A2)]
\item[(A1)] the field $f$ is almost surely continuous on $T$;

\item[(A2)] for some $\delta>0$ and $|s-t |<\delta$, the mean and variance
functions satisfies%
\[
\vert\sigma( t) -\sigma( s)
\vert
+\vert\mu( t) -\mu( s) \vert
\leq
\kappa_{H}|s-t|^{\beta};
\]

\item[(A3)] for some $\delta>0$ and $|s-s^{\prime}|<\delta$,
$|t-t^{\prime
}|<\delta$ the correlation function of $f$ satisfies%
\[
\vert r( t,s) -r( t^{\prime},s^{\prime
})
\vert\leq\kappa_{H}[|t-t^{\prime}|^{\beta}+|s-s^{\prime
}|^{\beta}];
\]

\item[(A4)] $0\in T$. There exist $\kappa_{0}$ and $\omega_{d}$ such
that, for
any $t\in T$ and $\varepsilon$ small enough,
\[
m\bigl(B(t,\varepsilon)\cap T\bigr)\geq\kappa_{0}\varepsilon^{d}\omega_{d},
\]
where $m$ is the Lebesgue measure, $B(t,\varepsilon)=\{s\dvtx
\vert
t-s\vert\leq\varepsilon\}$ and $\omega_{d}=m(
B(0,1)) $.
\end{longlist}

The assumption is that $0\in T$ is of no real consequence and is
adopted only for notational convenience.
\begin{theorem}
\label{ThmGeneral}Suppose that $f\dvtx T\to\real$ is a Gaussian
random field
satisfying conditions \textup{(A1)--(A4)} above. Then, algorithm~\ref{algorithm2}
% in Section~\ref{SecHolder}
provides a FPRAS for $w( b) $.
\end{theorem}

The polynomial rate of the intrinsic complexity bound inherent in this result
is discussed in Section
\ref{SecHolder}, along with similar rates in results to follow.
The conditions of the
previous theorem are weak and hold for virtually all applied settings involving
continuous Gaussian fields on compact sets.\vadjust{\goodbreak}

Not surprisingly, the complexity bounds of our algorithms can be improved
upon
under additional assumptions on $f$. For example, in the case of finite
fields (i.e., when $T$ is finite) with a nonsingular covariance
matrix, we
can show that the complexity of the algorithm is actually bounded as
$b\nearrow\infty$. We summarize this in the next result, whose proof,
which is
given in Section~\ref{SecFinite}, is useful for understanding the main ideas
behind the general procedure.\vspace*{-2pt}
\begin{theorem}
\label{ThmFinite}Suppose that $T$ is a finite set, and $f$ has a nonsingular
covariance matrix over $T\times T$. Then Algorithm~\ref{algorithm1}
%, in Section~\ref{SecFinite},
provides a FPRAS with $q=0$.\vspace*{-2pt}
\end{theorem}

As we indicated above, the strategy behind the discrete case provides the
basis for the general case. In the general situation, the underlying
idea is to discretize the field with an
appropriate sampling (discretization) rule that depends on the level
$b$ and the
continuity characteristics of the field. The number of sampling
points grows as $b$ increases, and the complexity of the algorithm is
controlled by finding a good sampling rule. There is a trade-off between
the number of points sampled, which has a direct impact on the
complexity of the algorithm, and the fidelity of the discrete
approximation to the continuous field. Naturally, in the presence of enough
smoothness and regularity, more information can be obtained with
the same sample size. This point is illustrated in the next
result, Theorem~\ref{ThmHom}, which considers smooth, homogeneous
fields. Note that in addition to controlling the
error induced by discretizing the field, \textit{the variance is strongly
controlled} and the \textit{discretization rule is optimal}, in a
sense explained in Section~\ref{SecSmooth}. For Theorem~\ref{ThmHom} we require the following additional regularity conditions:

\begin{longlist}[(B2)]
\item[(B1)] $f$ is homogeneous and almost surely twice continuously
differentiable;

%and $\partial_{i}C(0)=\partial_{ijk}C(0)=0$, for all
%$i,j,k=1,\ldots,d$.

\item[(B2)] $0 \in T\subset\real^{d}$ is a $d$-dimensional convex
set with nonempty
interior. Denoting its boundary by $\partial T$, assume that $\partial
T$ is a
$(d-1)$-dimensional manifold without boundary. For any $t\in T$, assume that
there exists $\kappa_{0}>0$ such that
\[
m\bigl(B(t,\varepsilon)\cap T\bigr)\geq\kappa_{0} \varepsilon^{d}
\]
for any $\varepsilon<1$, where $m$ is Lebesgue measure.\vspace*{-2pt}
\end{longlist}
\begin{theorem}
\label{ThmHom}Let $f$ satisfy conditions \textup{(B1)} and
\textup{(B2)}. Then Algorithm~\ref
{algorithm3}
% in Section~\ref{SecSmooth}
provides a FPRAS. Moreover, the underlying estimator is
strongly efficient and there exists a discretization scheme for $f$
which is
optimal in the sense of Theorem~\ref{PropOpt}.\vspace*{-2pt}%Section~\ref{SecSmooth}.
\end{theorem}

The results stated in Theorem~\ref{ThmHom} are stronger than those in
Theorem~\ref{ThmGeneral}. This is because conditions (B1) and (B2) are
much stronger than conditions~(A1)--(A4). The\vadjust{\goodbreak} structure present in Theorem
\ref{ThmHom} allows us to carry out a~more refined complexity
analysis. Using smoothness and homogeneity, the conditional
distribution of the random field given a high excursion can be
described quite precisely in an asymptotic sense using its derivatives.
In our analysis we take advantage of such a conditional description,
which is not available for H\"older continuous fields. On the other
hand, it might be possible that the algorithms developed for
Theorem~\ref{ThmGeneral}, or closely related variations, are in fact strongly
efficient for certain H\"older continuous fields. We leave this more
refined analysis to future study.

%s4 ###
\section{Importance sampling}\label{SecIS}

Importance sampling is based on the basic identity, for fixed
measurable $B$,
%
%e4.1 ###
%
\begin{equation} \label{Id1}%
P( B) =\int\indic( \omega\in B) \,dP(
\omega)
=\int\indic( \omega\in B) \,\frac{dP}{dQ}( \omega
)
\,dQ( \omega) ,
\end{equation}
where we assume that the probability measure $Q $ is such
that $Q( \cdot\cap B) $ is absolutely continuous with
respect to
the measure $P( \cdot\cap B) $. If we use $E^{Q} $ to denote
expectation under $Q $, then~(\ref{Id1}) trivially yields
that the random variable
\[
L( \omega) =\indic( \omega\in B) \,\frac
{dP}{dQ}(
\omega)
\]
is an unbiased estimator for $P( B) >0$ under the measure
$Q$, or, symbolically, $E^{Q}L=P( B) $.

An averaged importance sampling estimator based on the measure $Q$,
which is often referred as an \textit{importance sampling
distribution} or a \textit{change-of-measure}, is obtained by
simulating $n$
i.i.d. copies $L^{(1)},\ldots,L^{(n)}$ of $L$ under~$Q $ and
computing the empirical average $\widehat
{L}_{n}=(L^{(1)}+\cdots+L^{(n)})/n$. A
central\vspace*{1pt} issue is that of selecting $Q$ in
order to minimize the variance of $\widehat{L}_{n}$. It is easy to
verify that
if $\mathcal{Q}^{\ast}(\cdot)=P(\cdot|B)=P(\cdot\cap B)/P(B)$,
then the
corresponding estimator has zero variance. However, $\mathcal{Q}^{\ast
}$ is clearly a change of measure that is of no practical value, since
$P( B) $---the quantity that we are attempting to
evaluate in
the first place---is unknown. Nevertheless, when constructing a~good
importance sampling distribution for a family of sets $\{B_{b}\dvtx
b\geq b_{0}\}$
for which $0<P( B_{b}) \rightarrow0$ as $b\rightarrow
\infty$, it
is often useful to analyze the asymptotic behavior of $\mathcal
{Q}^{\ast}$ as $P( B_{b}) \rightarrow0$ in order to
guide the construction of a good~$Q$.

We now describe briefly how an efficient importance sampling estimator for
$P( B_{b}) $ can also be used to estimate a large class of
conditional expectations given~$B_{b}$. Suppose that a single
replication of
the corresponding importance sampling estimator,
\[
L_{b} \definedas\indic( \omega\in
B_{b}) \,dP/dQ
\]
can be generated in $O( {\log}\vert P(
B_{b}) \vert^{q_{0}})$ function evaluations,
for some $q_{0}>0$, and that
\[
\operatorname{Var}( L_{b}) =O({[P( B_{b})]^{2}\log}\vert
P( B_{b}) \vert^{q_{0}}).\vadjust{\goodbreak}
\]
These assumptions imply that by taking the average of i.i.d.
replications of~$L_{b}$ we obtain a FPRAS.

Fix $\beta\in(0,\infty)$, and let $\mathcal{X}( \beta
,q) $ be
the class of random variables $X$ satisfying $0\leq X\leq\beta$ with
%
%e4.2 ###
%
\begin{equation}
\label{cond11}
E[ X\vert B_{b}]=\Omega\lbrack1/\log(P( B_{b})
)^{q}].
\end{equation}
Then, by noting that%
%
%e4.3 ###
%
\begin{equation} \label{eq2}%
\frac{E^{Q}( XL_{b}) }{E^{Q}( L_{b})
}=E[
X\vert B_{b}]=\frac{E[X;B_{b}]}{P(B_{b}) },
\end{equation}
it follows easily that a FPRAS can be obtained by constructing the natural
estimator for $E[ X\vert B_{b}]$; that is, the ratio of the
corresponding averaged importance sampling estimators suggested by the ratio
in the left of~(\ref{eq2}). Of course, when $X$ is difficult to simulate
exactly, one must assume the bias $E[X;B_{b}]$ can be reduced in polynomial
time. The estimator is naturally biased, but the discussion on FPRAS on biased
estimators given in Section~\ref{SecSet} can be directly applied.

In the context of Gaussian random fields, we have that $B_{b}=\{f^{\ast
}>b\}$,
and one is very often interested in random variables $X$ of the form
$X=\Gamma
( f) $, where $\Gamma\dvtx C( T) \to\real
$, and
$C(T) $ denotes the space of continuous functions on~$T$. Endowing
$C(T) $ with the uniform topology, consider functions
$\Gamma$ that are nonnegative and bounded by a positive
constant. An archetypical example is given by
the volume of (conditioned) high-level excursion sets
with $\beta= m(T)$ is known to satisfy~(\ref{cond11}).
However, there are many other examples
of $\mathcal{X}( \beta,q) $ with
$\beta=m( T) $ which satisfy~(\ref{cond11}) for a suitable $q$,
depending on the regularity
properties of the field. In fact, if the mean and covariance properties of
$f $ are H\"{o}lder continuous, then, using similar
techniques as those given in the arguments of Section~\ref{SecHolder},
it is
not difficult to see that $q$ can be estimated.

In case that $\Gamma(f)$ is not bounded, the analysis is usually
case-by-case. In particular, we need to investigate
\[
E^Q(\Gamma^2(f)L_b^2)= E(\Gamma^2(f)L_b |B_b)P(B_b).
\]
We provide a brief calculation for the case of the conditional
overshoot, that is, $\Gamma(f) = f^* - b$ and $B_b = \{f^*>b\}$. We
admit the change of measure defined later in~(\ref{bkp}).
Then, given $\{f^*>b\}$, $\Gamma^2(f)$ and $L_b$ are negatively
correlated (the higher the overshoot is, the larger the excursion set
is), and we can obtain that
\[
E^Q(\Gamma^2(f)L_b^2)\leq E(\Gamma^2(f)|B_b)E(L_b).
\]
Conditional on the occurrence of $\{f^*>b\}$, $b\Gamma(f)$
asymptotically follows an exponential distribution. Therefore,
$E(\Gamma^2(f)|B_b) = (1+o(1))E^2(\Gamma(f)|B_b)$. Together with the
FPRAS of $L_b$ in computing $P(B_b)$, $\Gamma(f)L_b$ is an FPRAS to
compute the conditional overshoot.\vadjust{\goodbreak} The corresponding numerical examples
are given in Section~\ref{SecNum}. Two key steps involve the analyses
of the conditional correlation of $\Gamma^2 (f)$ and $L_b$ and the
conditional distribution of $\Gamma(f)$ given~$B_b$.

Thus, we have that a FPRAS based importance sampling algorithm for~$w(b)$
would typically also yield a polynomial time
algorithm for functional
characteristics of the conditional field given high level
excursions. Since this is a~very important and novel application, we
devote the remainder of this paper to
the development of
efficient importance sampling algorithms for~$w(b)$.

%s5 ###
\section{The basic strategy: Finite fields}\label{SecFinite}

In this section we develop our main ideas in the setting in which $T$ is
a finite set of the form $T=\{t_{1},\ldots,t_{M}\}$. To emphasize the discrete
nature of our algorithms in this section, we write $X_{i}=f(
t_{i}) $ for $1,\ldots,M$ and set ${X}=( X_{1},\ldots,X_{M}%
) $. This section is mainly of an expository nature, since much of
it has already appeared in~\cite{ABL08Gaussian}.
Nevertheless, it is included here as a useful guide to the intuition behind
the algorithms for the continuous case.

We have already noted that in order to design an efficient importance sampling
estimator for $w( b) =P( \max_{1\leq i\leq M}%
X_{i}>b) $ it is useful to study the asymptotic conditional
distribution of ${X}$, given that $\max_{1\leq i\leq M}X_{i}>b$. Thus,
we begin with some basic large deviation results.
\begin{proposition}
\label{PropFinite}%
For any set of random variables $X_1,\ldots,X_M$,
\[
\max_{1\leq i\leq M}P( X_{i}>b) \leq P\Bigl( \max
_{1\leq i\leq
M}X_{i}>b\Bigr) \leq\sum_{j=1}^{M}P( X_{j}>b) .
\]
Moreover, if the $X_j$ are mean zero, Gaussian, and the correlation
between~$X_{i}$ and $X_{j}$ is strictly less than 1, then%
\[
P( X_{i}>b,X_{j}>b) =o\bigl( \max
[P(X_{i}>b),P(X_{j}>b)]\bigr).
\]
Thus, if the associated covariance matrix of ${X}$ is
nonsingular,
\[
w( b) =\bigl(1+o( 1) \bigr)\sum_{j=1}^{M}P( X_{j}>b) .
\]
\end{proposition}
\begin{pf}
The lower bound in the first display is trivial,
and the upper bound follows by the union bound. The second display follows
easily by working with the joint density of a bivariate Gaussian distribution
(e.g.,~\cite{Berman93,LLR})
and the third claim is a direct consequence of the inclusion--exclusion
principle.
\end{pf}

As noted above, $\mathcal{Q}^{\ast}$ corresponds
to the conditional distribution of ${X}$ given that $X^{\ast}%
\triangleq\max_{1\leq i\leq M}X_{i}>b$. It follows from Proposition
\ref{PropFinite} that,
conditional on \mbox{$X^{\ast}>b$}, the probability
that two or more $X_j$ exceed $b$
is negligible. Moreover, it also follows that
\[
P( X_{i}=X^{\ast}\vert X^{\ast}>b) =\bigl(
1+o( 1) \bigr) \frac{P( X_{i}>b) }{\sum
_{j=1}%
^{M}P( X_{j}>b) }.
\]
The following corollary now follows as an easy consequence of these
observations.
\begin{corollary}%
\label{corollary-dtv}
\[
d_{\mathrm{TV}}(\mathcal{Q}^{\ast},\mathcal{Q})\rightarrow0
\]
as $b\rightarrow\infty$, where $d_{\mathrm{TV}}$ denotes the
total variation norm, and $\mathcal{Q} $ is defined,
for Borel $B\subset\real^{M}$, as
\[
\mathcal{Q}( {X}\in B) =\sum_{i=1}^{M}p( i,b)
P[ {X}\in B\vert X_{i}>b],
\]
where%
\[
p( i,b) =\frac{P( X_{i}>b) }{\sum
_{j=1}^{M}P(
X_{j}>b) }.
\]
\end{corollary}
\begin{pf}
Pick an arbitrary Borel $B$. Then we have that%
\begin{eqnarray*}
\mathcal{Q}^{\ast}( {X}\in B) &=& \frac{P[{X}\in
B,\max_{1\leq i\leq M}X_{i}>b]}{w( b) }\leq\sum
_{i=1}^{M}%
\frac{P[{X}\in B,X_{i}>b]}{w( b) }\\
&=& \sum_{i=1}^{M}P[{X}\in B|X_{i}>b]\frac{p( i,b)
}{(1+o( 1) )}.
\end{eqnarray*}
The above, which follows from the union bound and the last part
of Proposition~\ref{PropFinite} combined with the definition of
$\mathcal{Q} $, yields that for each $\varepsilon>0$ there
exists $b_{0}$ (independent of $B$) such that, for all $b\geq b_{0}$,%
\[
\mathcal{Q}^{\ast}( {X}\in B) \leq\mathcal{Q}(
{X}\in B) /(1-\varepsilon).
\]
The lower bound follows similarly, using the inclusion--exclusion
principle and the second part of Proposition~\ref{PropFinite}.
\end{pf}

Corollary~\ref{corollary-dtv}
provides support for choosing $\mathcal{Q}$ as an importance sampling
distribution. Of course, we still have
to verify that the corresponding algorithm is a FPRAS. The importance
sampling estimator induced by $\mathcal{Q} $ takes the
form%
%
%e5.1 ###
%
\begin{equation} \label{EstL1}%
L_{b}=\frac{dP}{dQ}=\frac{\sum_{j=1}^{M}P( X_{j}>b)
}{\sum
_{j=1}^{M}\indic( X_{j}>b) }.
\end{equation}
Note that under $\mathcal{Q} $ we have that
$X^{\ast}>b$ almost surely, so the denominator in the expression for $L_{b}$
is at least 1. Therefore, we have that%
\[
E^{\mathcal{Q}}L_{b}^{2}\leq\Biggl( \sum_{j=1}^{M}P(
X_{j}>b)
\Biggr) ^{2},
\]
and by virtue of Proposition~\ref{PropFinite} we conclude (using
$\operatorname{Var}_{\mathcal{Q}} $ to denote the variance under
$\mathcal{Q} $) that%
\[
\frac{\operatorname{Var}_{\mathcal{Q}}( L_{b}) }{P( X^{\ast
}>b)
^{2}}\rightarrow0
\]
as $b\rightarrow\infty$. In particular, it follows that $L_{b}$ is
strongly efficient.

Our proposed algorithm can now be summarized as follows.
\begin{algorithm}
%{\it
\label{algorithm1}
There are two steps in the algorithm:

\textsc{Step} (1).
Simulate $n$ i.i.d. copies ${X}^{( 1)
},\ldots,{X}^{( n) }$ of ${X}$ from the distribution~$\mathcal{Q} $.

\textsc{Step} (2).
Compute and output
\[
\widehat{L}_{n}=\frac{1}{n}\sum_{i=1}^{n}L_{b}^{( i) },
\]
where $L_{b}^{( i) }=\sum_{j=1}^{M}P(X_{j}^{(i)}>b)/\sum
_{j=1}^{M}\indic(X_{j}^{(i)}>b)$.
\end{algorithm}
%
%}

Since the generation of $L_{i}$ under $\mathcal{Q} $
takes $O( M^{3}) $ function evaluations we conclude, based
on the
analysis given in Section~\ref{SecSet}, that Algorithm \ref
{algorithm1} is a~FPRAS with $q=0$. This implies Theorem~\ref{ThmFinite}, as promised.

%s6 ###
\section{\texorpdfstring{A FPRAS for H\"older continuous Gaussian fields}{A FPRAS for Holder continuous Gaussian fields}}\label{SecHolder}
In this section we shall describe the algorithm and the analysis behind
Theorem~\ref{ThmGeneral}. Throughout the section, unless stated
otherwise, we
assume conditions (A1)--(A4) of Section~\ref{SecMain}.

There are two issues related to the complexity analysis. First, since
$f$ is
assumed continuous, the entire field cannot be generated in a (discrete)
computer, and so the algorithm used in the discrete case needs
adaptation. Once
this is done, we need to carry out an appropriate variance analysis.

Developing an estimator directly applicable to the continuous
field will be carried out in Section~\ref{SubCE}.
This construction will not only be
useful when studying the performance of a suitable discretization, but
will also help to explain some of the
features of our discrete construction. Then,
in Section~\ref{SubsecAAHolder}, we introduce a discretization approach
and study the bias caused by the discretization. In addition, we
provide bounds on the variance of this discrete importance sampling estimator.

%s6.1 ###
\subsection{A continuous estimator}\label{SubCE}

We start with a change of measure motivated by the discrete case in
Section~\ref{SecFinite}. A natural approach is to consider an importance
sampling strategy analogous to that of Algorithm~\ref{algorithm1}. The
continuous adaptation
involves first sampling $\tau_{b}$ according to the probability measure
%
%e6.1 ###
%
\begin{equation} \label{cont1}%
P( \tau_{b}\in\cdot) =\frac{E[m( A_{b}\cap\cdot
)
]}{E[m( A_{b}) ]},
\end{equation}
where $A_{b}=\{t\in T\dvtx f( t) >b\}$.
% Given that $f(
%measurable.
%The distribution of $\tau_{b}$ is well defined.
The idea of
introducing $\tau_{b}$ in the continuous setting is not necessarily to locate
the point at which the maximum is achieved, as was the situation
in the discrete case.
Rather, $\tau_{b}$ will be used
to find a random point which has a reasonable probability of being in
the excursion set $A_b$. (This probability will tend to be higher if
$f$ is nonhomogenous.) This relaxation will prove useful in
the analysis of the algorithm.
Note that $\tau_{b}$, with the distribution indicated in
(\ref{cont1}), has a density function (with respect to Lebesgue
measure) given
by
\[
h_{b}(t)=\frac{P(f(t)>b)}{E[m(A_{b})]},%
\]
and that we also can write
\[
E[m( A_{b}) ]=E\int_{T}\indic\bigl( f( t)
>b\bigr)\,
dt=\int_{T}P\bigl( f( t) >b\bigr) \,dt=m( T)
P\bigl( f( U) >b\bigr),
\]
where $U$ is uniformly distributed over $T$.

Once $\tau_{b}$ is generated, the natural continuous adaptation corresponding
to the strategy described by Algorithm~\ref{algorithm1} proceeds by
sampling $f$
conditional on $f(\tau_{b})>b$. Note that if we use $\bar{Q}$ to denote
the change-of-measure induced by such a continuous sampling strategy,
then the
corresponding importance sampling estimator takes the form%
\[
\bar{L}_{b}=\frac{dP}{d\bar{Q}}=\frac{E[m( A_{b})
]}{m(A_{b})}.
\]
The second moment of the estimator then satisfies
%
%e6.2 ###
%
\begin{eqnarray}
\label{variance-moment} \label{VarNHolder}%
E^{\bar{Q}}[(\bar{L}_{b})^{2}] &=& E( \bar{L}_{b};A_{b}\neq
\varnothing) \nonumber\\[-8pt]\\[-8pt]
&=& E[m(A_{b})]P( f^{\ast}>b) E[m(A_{b})^{-1}|A_{b}%
\neq\varnothing].\nonumber
\end{eqnarray}
Unfortunately, it is easy to construct examples for which $E^{\bar{Q}%
}[( \bar{L}_{b}) ^{2}]$ is infinite.
For instance, consider a homogeneous and twice differentiable random
field with zero mean and unit variance living on $T=[0,1]^d$. Using the
Slepian model, discussed in Section~\ref{SecSmooth}, it follows that
the asymptotic distribution of the overshoot given $\{f^* > b\}$
satisfies
\[
b(f^* - b) \rightarrow S,
\]
weakly as $b\rightarrow\infty$ where $S$ is an exponential random
variable. Consequently, the distribution of $m(A_b)$ given $m(A_b)>0$ satisfies
\[
m(A_b) \rightarrow\kappa b^{-d}S^{d/2}
\]
for some constant $\kappa$.
Therefore, the second moment in~(\ref{variance-moment}) is infinity as
long as $d\geq2$. This example suggests that the construction of the
change of measure needs to be modified slightly.

%Part of the reason for which one obtains a variance blow up in
%(\ref{VarNHolder}) in the continuous case is because we are using $
%to \textquotedblleft locate\textquotedblright the position of the
%excursion
%(by selecting a point in its interior) in a greedy way. In the
%multivariate
%Gaussian setting studied in Section~\ref{SecFinite}, the value of the
%index at
%which the maximum is achieved is very close to $\tau_{b}$ in a very
%strong
%sense (total variation) as $b\rightarrow\infty$. In the continuous
%setting, on
%the other hand, in order to obtain such a strong type of aproximation
%one
%would really need to consider the local correlation structure of the
%field.
%So, instead of trying to locate the maximum of the field using $
%which would involve considering local correlation (and even boundary)
%effects
%in the distribution of $\tau_{b}$, we shall simply take the point of
%view that
%$\tau_{b}$ is just attempting to locate the excursion set in a
%reasonable way;
%not necessarily insisting in placing a point in the excursion set
%precisely at
%$\tau_{b}$.

Extreme value theory considerations similar to those explained in the
previous paragraph give that the
overshoot of $f$ over a given level $b$ will be of order $\Theta(1/b)$.
Thus, in order to keep
$\tau_{b}$ reasonably close to the excursion set, we shall also
consider the
possibility of an undershoot of size $\Theta( 1/b) $
right at
$\tau_{b}$. As we shall see, this relaxation will allow us to
prevent the variance in~(\ref{variance-moment}) becoming infinite.
Thus, instead of~(\ref{cont1}), we shall consider
$\tau_{b-a/b}$ with density
%
%e6.3 ###
%
\begin{equation} \label{cont2}
h_{b-a/b}(t)=\frac{P(f(t)>b-a/b)}{E[m(A_{b-a/b})]}
\end{equation}
for some $a>0$. To ease on later notation, write
\[
\gamma_{a,b}\triangleq b-a/b,\qquad \tau_{\gamma_{a,b}}=\tau_{b-a/b}.
\]

Let $Q^{\prime}$ be the change of measure induced by sampling
$f$ as follows. Given~$\tau_{\gamma_{a,b}}$, sample $f(\tau
_{\gamma_{a,b}})$ conditional on $f(\tau_{\gamma_{a,b}})>\gamma
_{a,b}$. In
turn, the rest of~$f$ follows its conditional
distribution (under the nominal, or original, measure) given the
observed value
$f(\tau_{\gamma_{a,b}})$. We then have that the corresponding Radon--Nikodym
derivative is
%
%e6.4 ###
%
\begin{equation}\label{LR}%
\frac{dP}{dQ^{\prime}}=\frac{E[m(A_{\gamma_{a,b}})]}{m(A_{\gamma_{a,b}})},
\end{equation}
and the importance sampling estimator $L_{b}^{\prime}$ is
%
%e6.5 ###
%
\begin{equation} \label{bkp}%
L_{b}^{\prime}=\frac{dP}{dQ^{\prime}}\indic( A_{b}\neq
\varnothing)
=\frac{E[m(A_{\gamma_{a,b}})]}{m(A_{\gamma_{a,b}})}\indic\bigl(
m(A_{b})>0\bigr).
\end{equation}
Note that we have used the continuity of the field in order to write
$\{m( A_{b}) >0\}=\{A_{b}\neq\varnothing\}$ almost surely.
The motivation behind this choice lies in the fact
that since $m(A_{\gamma_{a,b}})>m( A_{b}) >0$, the
denominator may now be big enough to control the second moment of the
estimator. In particular, we consider the homogeneous and twice
differentiable field mentioned previously. Given $m(A_b)>0$,
$m(A_{\gamma_{a,b}})$ is asymptotically lower bounded by $\kappa
a^{d/2}b^{-d}$.
As we shall see, introducing the undershoot of size $a/b$ will be
very useful in the technical development both in the remainder of this section
and in Section~\ref{SecSmooth}. In addition, its introduction also provides
insight into the appropriate form of the estimator needed when
discretizing the field.

%s6.2 ###
\subsection{Algorithm and analysis}\label{SubsecAAHolder}

We still need to face the problem of generating $f$ in a computer. Thus we
now concentrate on a suitable discretization scheme, still having in
mind the
change of measure leading to~(\ref{LR}). Since our interest is to ultimately
design algorithms that are efficient for estimating expectations such as
$E[\Gamma( f) |f^{\ast}>b]$, where $\Gamma$
may be a functional of the whole field, we shall use a global
discretization scheme.

Consider ${U}=(U_{1},\ldots,U_{M})$ where $U_{i}$ are i.i.d. uniform
random variables taking values in $T$ and independent of the field $f$.
Set $T_{M}=\{U_{1},\ldots,U_{M}\}$ and
$X_{i}=X_{i}( U_{i}) =f( U_{i}) $ for $1\leq
i\leq
M$. Then\vspace*{-1pt} $X=(X_{1},\ldots,X_{m})$ (conditional on ${U}$) is a multivariate
Gaussian random vector with conditional means $\mu(U_{i})\definedas
E(X_{i}|U_{i})$ and covariances $\mathcal{C}(
U_{i},U_{j}) \definedas \Cov(
X_{i},X_{j}\vert U_{i},U_{j}) $. Our strategy is
to approximate $w(b)$ by
\[
w_{M}(\gamma_{a,b})=P\Bigl(\max_{t\in T_{M}}f(t)>\gamma_{a,b}\Bigr)=E\Bigl[P\Bigl(\max
_{1\leq
i\leq M}X_{i}>\gamma_{a,b}\big|{U}\Bigr)\Bigr].
\]
Given the development in Section~\ref{SecFinite}, it might not be surprising
that if we can ensure that $M=M( \varepsilon,b) $ is polynomial
in $1/\varepsilon$ and $b$, then we shall be in a good position to
develop a
FPRAS. The idea is to apply an importance sampling strategy similar to
that we
considered in the construction of $L_{b}^{\prime}$ of~(\ref{bkp}),
but this time it will be conditional on ${U}$. In view of our
earlier discussions, we propose sampling from $Q^{\prime\prime} $
defined via%
\[
Q^{\prime\prime}( \mbox{$ X$}\in B\vert
{U}) =\sum_{i=1}^{M}p_{{U}}( i)
P[ \mbox{$ X$}\in B\vert X_{i}>\gamma_{a,b},{U}%
],
\]
where
\[
p_{{U}}( i) =\frac{P( X_{i}>\gamma
_{a,b}\vert{U}) }{\sum_{j=1}^{M}P(
X_{j}>\gamma_{a,b}\vert{U}) }.
\]
We then obtain the (conditional) importance sampling estimator
%
%e6.6 ###
%
\begin{equation}\label{Est2}%
L_{b}({U})=\frac{\sum_{i=1}^{M}P( X_{i}>\gamma
_{a,b}\vert{U}) }{\sum_{i=1}^{M}\indic(X_{i}>\gamma
_{a,b})} \indic\Bigl(\max_{i=1}^M X_i > \gamma_{a,b}\Bigr).
\end{equation}
Note that\vspace*{1pt} the event $\{\max_{i=1}^M X_i > \gamma_{a,b}\}$ occurs with
probability 1 under $Q''$. Therefore, the indicator $I(\max_{i=1}^M
X_i > \gamma_{a,b})$ will be omitted when it does not cause confusion.
It is clear that
\[
w_{M}(\gamma_{a,b})=E^{Q^{\prime\prime}}[L_{b}({U})].
\]

Suppose for the moment that $M$, $a$ and the number of replications $n$ have
been chosen. Our future analysis will, in particular, guide the
selection of
these parameters. Then the procedure is summarized by the next algorithm.
\begin{algorithm}
\label{algorithm2}
%&\textbf{Algorithm 2}
The algorithm has three steps:

\textsc{Step} (1). Simulate ${U}^{( 1) },\ldots,{U}%
^{( n) }$ which are $n$ i.i.d. copies of the vector
${U}=(U_{1},\ldots,U_{M})$ described above.

\textsc{Step} (2). Conditional on each ${U}^{( i)
}$, for
$i=1,\ldots,n$, generate $L_{b}^{( i) }({U}^{(i)})$ as
described by~(\ref{Est2}) by considering the distribution of $\mbox{$
X$}^{( i) }({U}^{(i)})=(X_{1}^{( i)
}(U_{1}^{( i) }),\break\ldots, X_{M}^{( i)
}(U_{M}^{(
i) }))$. Generate\vspace*{1pt} the $\mbox{$ X$}^{( i)
}({U}^{(i)})$ independently so that at the end we obtain that the
$L_{b}^{( i) }({U}^{(i)})$ are $n$ i.i.d. copies of
$L_{b}({U})$.

\textsc{Step} (3). Output
\[
\widehat{L}_{n}\bigl({U}^{( 1) },\ldots,{U}%
^{( n) }\bigr)=\frac{1}{n}\sum_{i=1}^{n}L_{b}^{(
i)
}\bigl({U}^{(i)}\bigr).
\]
\end{algorithm}

%s6.3 ###
\subsection{\texorpdfstring{Running time of Algorithm \protect\ref{algorithm2}: Bias and variance control}{Running time of Algorithm 6.1: Bias and variance control}}\label{SubsecBiasVarianceHolder}

The remainder of this section is devoted to the analysis of the running
time of the
Algorithm~\ref{algorithm2}. The first step lies in estimating
the bias and second
moment of $L_{b}({U})$ under the change of measure induced by the
sampling strategy of the algorithm, which we denote by $Q^{\prime
\prime} $.
We start with a simple bound for the second moment.
\begin{proposition}\label{PropBasicBounds}
$\!\!\!$There exists a finite $\lambda_{0}$, depending
on $\mu_{T}={\max_{t\in T}}\vert\mu( t)
\vert
$ \textit{and} $\sigma_{T}^{2}=\max_{t\in T}\sigma^{2}(
t)
$, for which
\[
E^{Q^{\prime\prime}}[L_{b}({U})^{2}]\leq\lambda_{0}M^{2}P\Bigl(
\max_{t\in T}f( t) >b\Bigr) ^{2}.
\]
\end{proposition}
\begin{pf}
Observe that
\begin{eqnarray*}
&&
E^{Q^{\prime\prime}}[L_{b}({U})^{2}] \\
&&\qquad \leq E\Biggl( \Biggl(
\sum_{i=1}^{M}P( X_{i}>\gamma_{a,b}\vert
U_{i})
\Biggr) ^{2}\Biggr) \\
&&\qquad \leq E\Biggl( \Biggl( \sum_{i=1}^{M}\sup_{t\in T}P\bigl(
f(
U_{i}) >\gamma_{a,b}\vert U_{i}=t\bigr) \Biggr)
^{2}\Biggr) \\
&&\qquad = M^{2}\max_{t\in T}P\bigl( f( t) >\gamma_{a,b}
\bigr) ^{2}\\
&&\qquad \leq \lambda_{0}M^{2}\max_{t\in T}P\bigl( f( t)
>b\bigr)
^{2}\\
&&\qquad \leq \lambda_{0}M^{2}P\Bigl( \max_{t\in T}f( t)
>b\Bigr)
^{2},
\end{eqnarray*}
which completes the proof.
\end{pf}

Next we obtain a preliminary estimate of the bias.
\begin{proposition}
\label{PropBias}
{For each }$M\geq1$ we have
\begin{eqnarray*}
\vert w( b) -w_{M}( \gamma_{a,b})
\vert& \leq &
E\bigl[ \exp\bigl( -Mm( A_{\gamma_{a,b}}) /m(T)\bigr)
;A_{b}\cap T\neq\varnothing\bigr] \\[-2pt]
&&{} + P\Bigl( \max_{t\in T}f( t) >\gamma
_{a,b},\max_{t\in
T}f( t) \leq b\Bigr) .
\end{eqnarray*}
\end{proposition}
\begin{pf}
Note that%
\begin{eqnarray*}
\vert w( b) -w_{M}( \gamma
_{a,b})
\vert
&\leq& P\Bigl(\max_{t\in T_{M}}f( t) \leq\gamma_{a,b},\max
_{t\in
T}f( t) >b\Bigr)\\[-2pt]
&&{}+P\Bigl(\max_{t\in T_{M}}f( t) >\gamma
_{a,b},\max_{t\in T}f( t) \leq b\Bigr).
\end{eqnarray*}
The second term is easily bounded by
\[
P\Bigl(\max_{t\in T_{M}}f( t) >\gamma_{a,b},\max_{t\in
T}f(
t) \leq b\Bigr)\leq P\Bigl(\max_{t\in T}f( t) >\gamma_{a,b}%
,\max_{t\in T}f( t) \leq b\Bigr).
\]
The first term can be bounded as follows:
\begin{eqnarray*}
&& P\Bigl( \max_{t\in T_{M}}f( t) \leq\gamma_{a,b},\max
_{t\in
T}f( t) >b\Bigr) \\[-2pt]
&&\qquad \leq E\bigl[ \bigl( P[ f( U_{i}) \leq
\gamma
_{a,b}\vert f ]\bigr) ^{M}\indic(A_{b}\cap
T\neq\varnothing)\bigr] \\[-2pt]
&&\qquad \leq E\bigl[ \bigl( 1-m(A_{\gamma_{a,b}})/m( T
) \bigr)
^{M};A_{b}\cap T\neq\varnothing\bigr] \\[-2pt]
&&\qquad \leq E\bigl[ \exp\bigl( -Mm(A_{\gamma_{a,b}})/m(
T) \bigr)
;A_{b}\cap T\neq\varnothing\bigr] .
\end{eqnarray*}
This completes the proof.
\end{pf}

The previous proposition shows that controlling the \textit{relative}
bias of
$L_{b}({U})$ requires finding bounds for%
%
%e6.7 ###
%
\begin{equation}\label{EqBndM}%
E\bigl[ \exp\bigl(-Mm(A_{\gamma_{a,b}})/m(T)\bigr);A_{b}\cap T\neq\varnothing\bigr]
\end{equation}
and%
%
%e6.8 ###
%
\begin{equation}\label{EqBndDens}%
P\Bigl(\max_{t\in T}f( t) >\gamma_{a,b},\max_{t\in T}f(
t) \leq b\Bigr),
\end{equation}
and so we develop these. To
bound~(\ref{EqBndM}) we take advantage of the importance sampling
strategy based on $Q^{\prime} $ introduced earlier in
(\ref{LR}). Write
%
%e6.9 ###
%
\begin{eqnarray} \label{APP1}
&& E\bigl[ \exp\bigl( -Mm(A_{\gamma_{a,b}})/m( T)
\bigr)
;m(A_{b})>0\bigr]\nonumber\\[-9pt]\\[-9pt]
&&\qquad =E^{Q^{\prime}}\biggl( \frac{\exp( -Mm(A_{\gamma
_{a,b}})/m(
T) ) }{m(A_{\gamma_{a,b}})};m(A_{b})>0\biggr) Em(
A_{\gamma_{a,b}}) .\nonumber
\end{eqnarray}
Furthermore, note that for each $\alpha>0$ we have%
%
%e6.10 ###
%
\begin{eqnarray} \label{AP1a}\quad
&& E^{Q^{\prime}}\biggl( \frac{\exp( -Mm( A_{\gamma
_{a,b}})
/m( T) ) }{m( A_{\gamma_{a,b}}) }%
;m(A_{b})>0\biggr)\nonumber\\[-2pt]
&&\qquad \leq\alpha^{-1}\exp\bigl( -M\alpha/m(T)\bigr)
\\[-2pt]
&&\qquad\quad{} +E^{Q^{\prime}}\biggl(
\frac{\exp( -Mm(A_{\gamma_{a,b}})/m(T)) }{m(
A_{\gamma
_{a,b}}) };m(A_{\gamma_{a,b}})\leq\alpha;m(A_{b})>0\biggr) .
\nonumber
\end{eqnarray}
The next result, whose proof is given in Section \ref
{SecPropSimplify}, gives
a bound for the above expectation.
\begin{proposition}
\label{PropSimplify}Let $\beta$ be as in conditions \textup{(A2)}
and \textup{(A3)}. For
any $v>0$,
there exist constants $\kappa,\lambda_{2}\in(0,\infty)$
[independent of
$a\in(0,1)$ and $b$, but dependent on $v$] such that if we select
\[
\alpha^{-1}\geq\kappa^{d/\beta}(b/a)^{(2+v)2d/\beta},
\]
and define $W$ such that $P( W>x) =\exp( -x^{\beta
/d}) $ for $x\geq0$, then
%
%e6.11 ###
%
\begin{eqnarray}
\label{SmallArea}
&&
E^{Q^{\prime}}\biggl( \frac{\exp( -Mm(A_{\gamma_{a,b}%
})/m( T) ) }{m(A_{\gamma_{a,b}})};m(A_{\gamma
_{a,b}}%
)\leq\alpha;m(A_{b})>0\biggr)
\nonumber\\[-8pt]\\[-8pt]
&&\qquad \leq EW^{2}\frac{m( T) }%
{\lambda_{2}^{2d/\beta}M}\biggl( \frac{b}{a}\biggr) ^{4d/\beta}.
\nonumber
\end{eqnarray}
\end{proposition}

The following result gives us a useful upper bound on
(\ref{EqBndDens}). The proof is given in Section~\ref{SecPropDensity}.
\begin{proposition}
\label{PropDensity}
Assume that conditions \textup{(A2)} and \textup{(A3)} are in force.
For any $v>0$, let $\rho=2d/\beta+dv+1$, where $d$ is the
dimension of $T$. There exist constants $b_{0},\lambda\in(
0,\infty) $ [independent of $a$ but depending on $\mu_{T}={\max
_{t\in
T}}\vert\mu( t) \vert$, $\sigma
_{T}^{2}=\max_{t\in
T}\sigma^{2}( t) $, $v$, the H\"older parameters $\beta$ and
$\kappa_{H}$] so that for all $b\geq b_{0}\geq1$ we have
%
%e6.12 ###
%
\begin{equation}\label{PropPart1}%
P\Bigl( \max_{t\in T}f(t)\leq b+a/b \big| \max_{t\in T}f(t)>b\Bigr)
\leq\lambda ab^{\rho}.
\end{equation}
Consequently,%
\begin{eqnarray*}
&& P\Bigl( \max_{t\in T}f( t) >\gamma_{a,b},\max_{t\in
T}f( t) \leq b\Bigr) \\
&&\qquad =P\Bigl( \max_{t\in T}f( t) \leq
b\big|\max_{t\in T}f(
t) >\gamma_{a,b}\Bigr) P\Bigl( \max_{t\in T}f(t)>\gamma
_{a,b}\Bigr)
\\
&&\qquad \leq\lambda ab^{\rho}P\Bigl( \max_{t\in
T}f(t)>\gamma_{a,b}\Bigr) .
\end{eqnarray*}
Moreover,
\[
P\Bigl( \max_{t\in T}f(t)>\gamma_{a,b}\Bigr) (1-\lambda ab^{\rho
})\leq
P\Bigl( \max_{t\in T}f(t)>b\Bigr) .
\]
\end{proposition}

Propositions~\ref{PropSimplify} and~\ref{PropDensity} allow us to prove
Theorem~\ref{ThmGeneral}, which is rephrased in the
form of the following theorem, which contains the detailed rate of complexity
and so the main result of this section.
\begin{theorem}
\label{PropFinalHolder} Suppose $f$ is a Gaussian random field satisfying
conditions \textup{(A1)--(A4)} in Section~\ref{SecMain}. Given any $v>0$, put
$a=\varepsilon
/(4\lambda b^{\rho})$ (where $\lambda$ and $\rho$ as in Proposition
\ref{PropDensity}), and\vadjust{\goodbreak} $\alpha^{-1}=\kappa^{d/\beta
}(b/a)^{(2+v)d/\beta
}$. Then, there exist $c,\varepsilon_{0}>0$ such that for all
$\varepsilon\leq\varepsilon_{0}$,
%
%e6.13 ###
%
\begin{equation}
\label{bias}\vert w( b) -w_{M}( \gamma
_{a,b})
\vert\leq w( b) \varepsilon,
\end{equation}
if $M=\lceil c\varepsilon^{-1}(b/a)^{(4+4v)d/\beta}
\rceil
$. Consequently, by our discussion in Section~\ref{SecSet} and the
bound on the second moment given in Proposition~\ref{PropBasicBounds}, it
follows that Algorithm~\ref{algorithm2} provides a FPRAS with running time
$O( (
M) ^{3}\times( M) ^{2}\times\varepsilon^{-2}
\delta^{-1}) $.
\end{theorem}
\begin{pf}
Combining~(\ref{PropBias}),~(\ref{APP1}) and
(\ref{AP1a}) with Propositions~\ref{PropBias}--\ref{PropDensity} we
have that
\begin{eqnarray*}
&& |w( b) -w_{M}( \gamma_{a,b}) |\\
&&\qquad \leq\alpha^{-1}\exp\bigl( -M\alpha/m(
T) \bigr) E[
m( A_{\gamma_{a,b}}) ]
\\
&&\qquad\quad{} +E[W^{2}]\frac{m( T)
}{\lambda_{2}^{2d/\beta}M}\biggl( \frac{b}{a}\biggr) ^{4d/\beta
}E[
m( A_{\gamma_{a,b}}) ] +\biggl( \frac{\lambda
ab^{\rho}%
}{1-\lambda ab^{\rho}}\biggr) w( b) .
\end{eqnarray*}
Furthermore, there exists a constant $K\in(0,\infty)$ such that%
\[
Em( A_{\gamma_{a,b}}) \leq K\max_{t\in T}P\bigl( f(
t) >b\bigr) m( T) \leq Kw( b) m(
T).
\]
Therefore, we have that%
\begin{eqnarray*}
\frac{|w( b) -w_{M}( \gamma_{a,b})
|}{w( b) }%
&\leq&\alpha^{-1}Km(T)\exp\bigl( -M\alpha/m( T)
\bigr)
\\
&&{} +E[W^{2}]K\frac{m( T) ^{2}}{\lambda
_{2}^{2d/\beta}M}\biggl( \frac
{b}{a}\biggr) ^{4d/\beta}+\biggl( \frac{\lambda ab^{\rho
}}{1-\lambda
ab^{\rho}}\biggr) .
\end{eqnarray*}
Moreover, since $a=\varepsilon/(4\lambda b^{\rho})$, we obtain that, for
$\varepsilon\leq1/2$,%
\begin{eqnarray*}
\frac{|w( b) -w_{M}( \gamma_{a,b})
|}{w( b) }%
&\leq&\alpha^{-1}Km(T)\exp\bigl( -M\alpha/m( T)
\bigr)
\\ &&{} +[EW^{2}]K\frac{m( T) ^{2}}{\lambda
_{2}^{2d/\beta}M}\biggl( \frac
{b}{a}\biggr) ^{4d/\beta}+\varepsilon/2.
\end{eqnarray*}
From the selection of $\alpha,M$ and $\theta$ it follows easily that
the first
two terms on the right-hand side of the previous display can be made
less than
$\varepsilon/2$ for all $\varepsilon\leq\varepsilon_{0}$ by taking
$\varepsilon_{0}$ sufficiently small.

The complexity count given in the theorem now
corresponds to the following estimates. The factor $O((
M) ^{3})$ represents the cost of a Cholesky factorization
required to
generate a single replication of a finite field of dimension $M$. In addition,
the second part of Proposition~\ref{PropBasicBounds} gives us that
$O(M^{2}\varepsilon^{-2}\delta^{-1})$ replications are required to control
the relative variance of the estimator.
\end{pf}

We now proceed to prove Propositions~\ref{PropSimplify} and~\ref{PropDensity}.

%s6.4 ###
\subsection{\texorpdfstring{Proof of Proposition \protect\ref{PropSimplify}}{Proof of Proposition 6.4}}
\label{SecPropSimplify}

We concentrate on the analysis of the left-hand side of (\ref
{SmallArea}). An
important observation is that conditional on the random variable $\tau
_{\gamma_{a,b}}$ with distribution%
\[
Q^{\prime}(\tau_{\gamma_{a,b}}\in\cdot)=\frac{E[
m(A_{\gamma_{a,b}}%
\cap\cdot)] }{E[ m(A_{\gamma_{a,b}})] }
\]
and, given $f(\tau_{\gamma_{a,b}})$, the rest of the field, namely
$(f(
t) \dvtx t\in T\setminus\{\tau_{\gamma_{a,b}}\})$ is
another Gaussian field
with a computable mean and covariance structure. The second term in
(\ref{AP1a}) indicates that we must estimate the probability that
$m(A_{\gamma
_{a,b}})$ takes small values under $Q^{\prime}$. For this purpose, we shall
develop an upper bound for%
%
%e6.14 ###
%
\begin{equation} \label{bndbasic}%
P\bigl( m(A_{\gamma_{a,b}})<y^{-1},m(A_{b})>0\vert
f(
t) =\gamma_{a,b}+z/\gamma_{a,b}\bigr)
\end{equation}
for $y$ large enough. Our arguments proceeds in two steps. For the
first, in order to study~(\ref{bndbasic}), we shall estimate the
conditional mean covariance of $\{f( s) \dvtx s\in T\}$,
given that $f(
t) =\gamma_{a,b}+z/\gamma_{a,b}$. Then, we use the fact that the
conditional field is also Gaussian and take advantage of general
results from
the theory of Gaussian random fields to obtain a bound for~(\ref{bndbasic}).
For this purpose we recall some useful results from the theory of Gaussian
random fields. The first result is due to Dudley~\cite{DUD73}.\vspace*{-2pt}
%
% and its proof can be
%consulted in~\cite{AdlerTaylor07}.
%
\begin{theorem}
\label{ThmA1} Let $\mathcal{U}$ be a compact subset of $\real^{n}$,
and let
$\{f_0( t) \dvtx t\in\mathcal{U}\}$ be a mean zero, continuous
Gaussian random field. Define the
canonical metric $d $ on $\mathcal{U}$ as%
\[
d( s,t) =\sqrt{E[f_0( t) -f_0( s)
]^{2}}
\]
and put $\operatorname{diam}( \mathcal{U}) =\sup_{s,t\in\mathcal
{U}}d(
s,t)$, which is assumed to be finite. Then there exists a finite
universal constant $\kappa>0$ such
that%
\[
E\Bigl[\max_{t\in\mathcal{U}}f_{0}( t) \Bigr]\leq\kappa\int_{0}
^{\mathrm{diam}( \mathcal{U}) /2}[\log
( \mathcal{N}(
\varepsilon) ) ]^{1/2}\,d\varepsilon,
\]
where the entropy $\mathcal{N}( \varepsilon) $
is the smallest number of $d$-balls of radius $\varepsilon$ whose
union covers $\mathcal{U}$.\vspace*{-2pt}
\end{theorem}

The second general result that we shall need is the so-called B--TIS
(Borel--Tsirelson--Ibragimov--Sudakov)
inequality~\cite{AdlerTaylor07,BOR,CIS}.\vspace*{-2pt}
\begin{theorem}
\label{ThmA2}Under the setting described in Theorem~\ref{ThmA1},%
\[
P\Bigl(\max_{t\in\mathcal{U}}f_0( t) -E\Bigl[\max_{t\in
\mathcal{U}}%
f_0( t) \Bigr]\geq b\Bigr)\leq\exp\bigl( -b^{2}/(2\sigma
_{\mathcal{U}%
}^{2})\bigr) ,
\]
where
\[
\sigma_{\mathcal{U}}^{2}=\max_{t\in\mathcal{U}}E[f_0^2( t)].\vspace*{-2pt}
\]
\end{theorem}

We can now proceed with the main proof. We shall
assume from now on that $\tau_{\gamma_{a,b}}=0$, since,\vadjust{\goodbreak} as will be
obvious from
what follows, all estimates hold uniformly over $\tau_{\gamma_{a,b}}
\in T$.
This is\vspace*{1pt} a consequence of the uniform H\"older assumptions (A2) and~(A3).
Define a new process $\widetilde f$
\[
\bigl(\widetilde
{f}( t) \dvtx t\in T\bigr)\equalinlaw\bigl( f( t
) \dvtx t\in
T\vert f( 0) =\gamma_{a,b}+z/\gamma_{a,b}\bigr).
\]
Note that we
can always write $\widetilde{f}( t) =\widetilde{\mu
}(
t) +g( t) $, where $g $ is a mean
zero Gaussian random field on $T$. We have that%
\[
\widetilde{\mu}( t) =E\widetilde{f}( t)
=\mu(
t) +\sigma( 0) ^{-2}C( 0,t) \bigl(
\gamma_{a,b}+z/\gamma_{a,b}-\mu( 0) \bigr) ,
\]
and that the covariance function of $\widetilde{f} $ is
given by
\[
C_{g}( s,t) =\Cov( g( s) ,g(
t)
) =C( s,t) -\sigma( 0) ^{-2}C(
0,s) C( 0,t) .
\]

The following lemma describes the behavior of $\widetilde{\mu}(
t) $ and $C_{g}( t,s) $.
\begin{lemma}
\label{LemSimple1}Assume that $|s|$ and $|t|$ small enough. Then the
following three conclusions hold:
\begin{longlist}
\item
There exist constants $\lambda_{0}$ and $\lambda_{1}>0$ such
that
\[
\vert\widetilde{\mu}( t) -(\gamma_{a,b}+z/\gamma
_{a,b})\vert\leq\lambda_{0}\vert t\vert^{\beta
}+\lambda
_{1}\vert t\vert^{\beta}(\gamma_{a,b}+z/\gamma_{a,b}),
\]
and for all $z\in(0,1)$ and $\gamma_{a,b}$ large enough,%
\[
|\widetilde{\mu}( s) -\widetilde{\mu}( t)
|\leq\kappa_{H}\gamma_{a,b}|s-t|^{\beta}.
\]

\item
\begin{eqnarray*}
\\[-26pt]
C_{g}( s,t) &\leq&2\kappa_{H}\sigma( t)
\sigma(
s) \{\vert t\vert^{\beta}+\vert s
\vert^{\beta
}+\vert t-s\vert^{\beta}\}.
\end{eqnarray*}

\item
\begin{eqnarray*}
\\[-26pt]
D_{g}( s,t) &=&\sqrt{E\bigl([ g( t) -g(
s) ] ^{2}\bigr)}\leq\lambda_{1}^{1/2}\vert t-s
\vert
^{\beta/2}.
\end{eqnarray*}
\end{longlist}
\end{lemma}
\begin{pf}
All three consequences follow from simple algebraic manipulations. The
details are omitted.
\end{pf}
\begin{proposition}
\label{PropTogether} For any $v>0$, there exist $\kappa$ and
$\lambda_{2}$, such that for all $t\in T$, $y^{-\beta/d}\leq\frac
{a^{2+v}%
}{\kappa b^{(2+v)}}$, $a$ sufficiently small, and $z>0$,
\[
P\bigl( m(A_{\gamma_{a,b}})^{-1}>y,m(A_{b})>0 | f
( t)
=\gamma_{a,b}+z/\gamma_{a,b}\bigr) \leq\exp( -\lambda
_{2}a^{2}%
y^{\beta/d}/b^{2}) .
\]
\end{proposition}
\begin{pf}
For notational simplicity, and without loss of generality,
we assume that $t=0$. First consider the
case that $z\geq1$. Then there exist $c_{1},c_{2}$ such that for all
$c_{2}y^{-\beta/d}<b^{-2-v}$ and $z>1$,%
\begin{eqnarray*}
&& P\bigl( m( A_{\gamma_{a,b}}\cap T)
^{-1}>y,m(A_{b})>0\big|f( 0) =\gamma_{a,b}+z/\gamma
_{a,b}\bigr) \\
&&\qquad \leq P\Bigl( \inf_{|t|<c_{1}y^{-1/d}}f(t)\leq\gamma
_{a,b}\big|f(
0) =\gamma_{a,b}+z/\gamma_{a,b}\Bigr) \\
&&\qquad =P\Bigl( \inf_{|t|<c_{1}y^{-1/d}}\widetilde{\mu
}( t)
+g(t)\leq\gamma_{a,b}\Bigr) \\[-2pt]
&&\qquad \leq P\biggl( \inf_{|t|<c_{1}y^{-1/d}}g(t)\leq-\frac
{1}{2\gamma_{a,b}%
}\biggr) .
\end{eqnarray*}
Now\vspace*{1pt} apply (iii) from Lemma~\ref{LemSimple1}, from which it follows
that $\mathcal{N}( \varepsilon) \leq c_{3}%
m(T)/\break\varepsilon^{2d/\beta}$ for some constant $c_{3}$. By Theorem
\ref{ThmA1}, $E(\sup_{|t|<c_{1}y^{-1/d}}f(t))=\break O(y^{-\beta/(2d)}\log
y)$. By
Theorem~\ref{ThmA2}, for some constant $c_{4}$,%
\begin{eqnarray*}
&& P\bigl( m( A_{\gamma_{a,b}}\cap T)
^{-1}>y,m(A_{b})>0|f( 0) =\gamma_{a,b}+z/\gamma
_{a,b}\bigr) \\[-2pt]
&&\qquad \leq P\biggl( \inf_{|t|<c_{1}y^{-1/d}}g(t)\leq-\frac
{1}{2\gamma_{a,b}%
}\biggr) \\[-2pt]
&&\qquad \leq\exp\biggl( -\frac{1}{c_{4}\gamma
_{a,b}^{2}y^{-\beta/d}}\biggr)
\end{eqnarray*}
for $c_{2}y^{-\beta/d}<b^{-2-v}$ and $z>1$.

Now consider the case
$z\in(0,1)$. Let $t^{\ast}$ be the global maximum of $f(t)$. Then,
\begin{eqnarray*}
&& P\bigl( m( A_{\gamma_{a,b}}\cap T) ^{-1}>y,m(A_{b}%
)>0|f( 0) =\gamma_{a,b}+z/\gamma_{a,b}\bigr) \\[-2pt]
&&\qquad \leq P\Bigl( \inf_{|t-t^{\ast
}|<c_{1}y^{-1/d}}f(t)<\gamma_{a,b},f(t^{\ast
})>b\big|f( 0) =\gamma_{a,b}+z/\gamma_{a,b}\Bigr) \\[-2pt]
&&\qquad \leq P\Bigl( \sup_{|s-t|<c_{1}y^{-1/d}}|f(s)-f(t)|>
a/b\big|f(
0) =\gamma_{a,b}+z/\gamma_{a,b}\Bigr) .
\end{eqnarray*}
Consider the new field $\xi(s,t)=g(s)-g(t)$ with parameter space
$T\times T$. Note that
\[
\sqrt{\operatorname{Var}(\xi(s,t))}=D_{g}(s,t)\leq\lambda_{1}|s-t|^{\beta/2}.
\]
Via basic algebra, it is not hard to show that the entropy of $\xi
(s,t)$ is bounded by
$\mathcal{N}_{\xi}( \varepsilon) \leq c_{3}m(T\times
T)/\varepsilon^{2d/\beta}$. In addition, from (i) of Lemma \ref
{LemSimple1}, we
have
\[
|\widetilde{\mu}( s) -\widetilde{\mu}( t)
|\leq\kappa_{H}\gamma_{a,b}|s-t|^{\beta}.
\]
Similarly, for some $\kappa>0$ and all $y^{-\beta/d}\leq\frac
{1}{\kappa
}( \frac{a}{b}) ^{2+v}$, $a<1$, there exists $c_{5}$ such that
\begin{eqnarray*}
&& P\bigl( m( A_{\gamma_{a,b}}\cap T) ^{-1}>y,m(A_{b}%
)>0|f( 0) =\gamma_{a,b}+z/\gamma_{a,b}\bigr) \\[-2pt]
&&\qquad
\leq P\Bigl( \sup_{|s-t|<c_{1}y^{-1/d}}|f(s)-f(t)|> a/b\big|f(
0) =\gamma_{a,b}+z/\gamma_{a,b}\Bigr) \\[-2pt]
&&\qquad =P\biggl( \sup_{|s-t|<c_{1}y^{-1/d}}|\xi(s,t)|>\frac
{a}{2b}\biggr) \\[-2pt]
&&\qquad \leq\exp\biggl( -\frac{a^{2}}{c_{5}b^{2}y^{-\beta
/d}}\biggr) .
\end{eqnarray*}
Combining the two cases $z>1$ and $z\in(0,1)$ and choosing $c_{5}$
large enough we have%
\begin{eqnarray*}
&& P\bigl( m( A_{\gamma_{a,b}}\cap T) ^{-1}>y,m(A_{b})>0
|
f( 0) =\gamma_{a,b}+z/\gamma_{a,b}\bigr) \\
&&\qquad \leq\exp\biggl( -\frac{a^{2}}{c_{5}b^{2}y^{-\beta
/d}}\biggr)
\end{eqnarray*}
for $a$ small enough and $y^{-\beta/d}\leq\frac{1}{\kappa}(
\frac{a}%
{b}) ^{2+v}$. Renaming the constants completes the proof.\vspace*{-2pt}
\end{pf}

The final ingredient needed for the proof of Proposition~\ref{PropSimplify}
is the
following lemma involving stochastic domination. The proof follows an
elementary argument and is therefore omitted.\vspace*{-2pt}
\begin{lemma}
\label{LemDominance}Let $v_{1} $ and $v_{2}$ be finite measures on
$\real$ and define $\eta_{j}(x)=\int
_{x}^{\infty}v_{j}( ds) $. Suppose that $\eta_{1}(
x) \geq\eta_{2}( x) $ for each $x\geq x_{0}$. Let
$(
h( x) \dvtx x\geq x_{0}) $ be a nondecreasing,
positive and
bounded function. Then,
\[
\int_{x_{0}}^{\infty}h( s) v_{1}(dx)\geq\int
_{x_{0}}^{\infty
}h( s) v_{2}( ds) .\vspace*{-2pt}
\]
\end{lemma}
\begin{pf*}{Proof of Proposition~\ref{PropSimplify}}
Note that
\begin{eqnarray*}
&& E^{Q^{\prime}}\biggl( \frac{\exp( -Mm( A_{\gamma
_{a,b}})
/m( T) ) }{m( A_{\gamma_{a,b}})
};m(
A_{\gamma_{a,b}}) \in(0,\alpha);m(A_{b})>0\biggr) \\[-2pt]
&&\qquad =\int_{T}\int_{z=0}^{\infty}E\biggl( \frac{\exp( -Mm(
A_{\gamma_{a,b}}) /m( T) ) }{m(
A_{\gamma_{a,b}}) };\\[-2pt]
&&\qquad\quad\hspace*{52pt}m( A_{\gamma_{a,b}}) \in
(0,\alpha);m(A_{b})>0\Big|
f( t) =\gamma_{a,b}+\frac
{z}{\gamma_{a,b}}\biggr)\\[-2pt]
&&\qquad\quad\hspace*{32pt}{}\times
P( \tau_{\gamma_{a,b}}\in dt) P\bigl(
\gamma_{a,b}[f( t) -\gamma_{a,b}]\in dz\vert
f(
t) >\gamma_{a,b}\bigr)\\[-2pt]
&&\qquad \leq\sup_{z>0,t\in T}E\biggl( \frac{\exp( -Mm(
A_{\gamma_{a,b}}) /m( T) ) }{m(
A_{\gamma_{a,b}}) };\\[-2pt]
&&\qquad\quad\hspace*{51pt}m( A_{\gamma_{a,b}}) \in
(0,\alpha);m(A_{b})>0\Big|
f( t) =\gamma_{a,b}+\frac
{z}{\gamma _{a,b}}\biggr) .
\end{eqnarray*}

Now define $Y( b/a) =(b/a)^{2d/\beta}\lambda
_{2}^{-d/\beta}W$ with
$\lambda_{2}$ as chosen in Proposition~\ref{PropTogether} and $W$
with distribution given by $P(W>x)=\exp(-x^{\beta/d})$. By Lemma \ref
{LemDominance},
\begin{eqnarray*}
&& \sup_{t\in T}E\biggl( \frac{\exp( -Mm( A_{\gamma_{a,b}%
}) /m( T) ) }{m( A_{\gamma
_{a,b}})
};\\[-2pt]
&&\qquad\hspace*{12pt} m( A_{\gamma_{a,b}}) \in(0,\alpha);m(A_{b})>0
\Big|
f( t) =\gamma_{a,b}+\frac{z}{\gamma_{a,b}}\biggr) \\[-2pt]
&&\qquad \leq E\bigl[ Y(b/a)\exp
\bigl(-MY(b/a)^{-1}/m(T)\bigr);Y(b/a)>\alpha^{-1}\bigr] .
\end{eqnarray*}
Now let $Z$ be exponentially distributed with mean 1 and independent of
$Y(b/a)$. Then we have (using the definition of the tail distribution
of $Z$
and Chebyshev's inequality)%
\begin{eqnarray*}
\exp\bigl( -MY(b/a)^{-1}/m( T) \bigr) &=& P\bigl(
Z>MY(b/a)^{-1}/m( T) \vert Y( b/a)
\bigr)
\\ &\leq&\frac{m( T) Y( b/a) }{M}.
\end{eqnarray*}
Therefore,
\begin{eqnarray*}
&& E\bigl[\exp\bigl( -MY(b/a)^{-1}/m( T) \bigr)
Y(b/a);Y(b/a)\geq\alpha^{-1}\bigr]
\\
&&\qquad
\leq\frac{m( T) }{M}E[Y(b/a)^{2};Y(b/a)\geq\alpha
^{-1}]\\
&&\qquad \leq\frac{m( T) }{M\lambda_{2}^{2d/\beta
}}\biggl( \frac{b}%
{a}\biggr) ^{4d/\beta}E(W^{2}),
\end{eqnarray*}
which completes the proof.
\end{pf*}
%

%s6.5 ###
\subsection{\texorpdfstring{Proof of Proposition \protect\ref{PropDensity}}{Proof of Proposition 6.5}}
\label{SecPropDensity}

We start with the following
result of Tsirel\-son~\cite{TsirelsonMaxDens}.

\begin{theorem}
\label{tsirelsontheorem} Let $f$ be a continuous separable Gaussian process
on a compact (in the canonical metric) domain $T$. Suppose that
$\operatorname{Var}(
f ) =\sigma$ is continuous
and that $\sigma( t) >0$ for $t\in T$. Moreover, assume that
$\mu=Ef$ is also continuous and $\mu(t) \geq0$ for all $t\in T$. Define
\[
\sigma_{T}^{2}\stackrel{\Delta}{=}\max_{t\in T}\operatorname{Var}(f(t))
\]
and set $F(x)=P\{\max_{t\in T}f(t)\leq x\}$. Then, $F $
is continuously differentiable on $\real$. Furthermore, let $y$ be
such that
$F(y)>1/2$, and define $y_{\ast}$ by
\[
F(y)=\Phi(y_{\ast}).
\]
Then, for all
$x>y$,
\[
F^{\prime}(x)\leq\Psi\biggl( \frac{xy_{\ast}}{y}\biggr) \biggl(
\frac{xy_{\ast}}{y}(1+2\alpha)+1\biggr) (1+\alpha),
\]
where
\[
\alpha=\frac{y^{2}}{x(x-y)y_{\ast}^{2}}.
\]
\end{theorem}

We can now prove the following lemma.
\begin{lemma}
\label{LemDens1}There exists a constant $A\in( 0,\infty
) $
independent of $a$ and $b\geq0$ such that%
%
%e6.15 ###
%
\begin{equation} \label{LemPart2}
P\Bigl( \sup_{t\in T}f(t)\leq b+a/b \big| \sup_{t\in T}f(t)>b\Bigr)
\leq aA
\frac{P( \sup_{t\in T}f( t) \geq b-1/b)
}{P(
\sup_{t\in T}f( t) \geq b) }.\hspace*{-35pt}
\end{equation}
\end{lemma}
\begin{pf}
By subtracting $\inf_{t\in
T}\mu(
t) >-\infty$ and redefining the level $b$ to be $b-\inf_{t\in
T}\mu(
t) $ we may simply assume that $Ef( t) \geq0$ so
that we
can apply Theorem~\ref{tsirelsontheorem}. Adopting the notation of Theorem
\ref{tsirelsontheorem}, first we pick~$b_{0}$ large enough so that
$F(
b_{0}) >1/2$ and assume that $b\geq b_{0}+1$. Now, let
$y=b-1/b$ and $F(y)=\Phi(y_{\ast})$. Note that there exists $\delta
_{0}%
\in(0,\infty)$ such that $\delta_{0}b\leq y_{\ast}\leq\delta
_{0}^{-1}b$ for all
$b\geq b_{0}$. This follows easily from the fact that
\[
\log{{\mathbb{P}}}\Bigl\{ \sup_{t\in T}f(t) > x\Bigr\} \sim
\log
\sup_{t\in T}{{\mathbb{P}}}\{ f(t) > x\} \sim-\frac
{x^{2}%
}{2\sigma_{T}^{2}}.
\]
On
the other hand, by Theorem~\ref{tsirelsontheorem} $F $ is
continuously differentiable, and~so%
%
%e6.16 ###
%
\begin{equation}\label{eqref}\quad
P\Bigl\{ \sup_{t\in T}f(t)<b+a/b \big| \sup_{t\in T}f(t)>b\Bigr\}
=\frac
{\int_{b}^{b+a/b}F^{\prime}(x) \,dx}{P\{ \sup_{t\in
T}f(t)>b\} }.
\end{equation}
Moreover,%
\begin{eqnarray*}
F^{\prime}(x) & \leq &\biggl( 1-\Phi\biggl( \frac{xy_{\ast
}}{y}\biggr)
\biggr) \biggl( \frac{xy_{\ast}}{y}\bigl(1+2\alpha( x)
\bigr)+1\biggr)
\frac{y_{\ast}}{y}\bigl(1+\alpha( x) \bigr)\\
& \leq &\bigl( 1-\Phi( y_{\ast}) \bigr) \biggl( \frac
{xy_{\ast
}}{y}\bigl(1+2\alpha( x) \bigr)+1\biggr) \frac{y_{\ast
}}{y}\bigl(1+\alpha
( x) \bigr)\\
& = &P\Bigl( \max_{t\in T}f( t) >b-1/b\Bigr) \biggl(
\frac{xy_{\ast}}{y}\bigl(1+2\alpha( x) \bigr)+1\biggr) \frac
{y_{\ast}}%
{y}\bigl(1+\alpha( x) \bigr).
\end{eqnarray*}
Therefore,
\begin{eqnarray*}
&&\int_{b}^{b+a/b}F^{\prime}(x) \,dx \\
&&\qquad \leq P\Bigl( \sup_{t\in T}f( t)
>b-1/b\Bigr) \int_{b}^{b+a/b}\biggl( \frac{xy_{\ast}}{y}\bigl(1+2\alpha
(
x) \bigr)+1\biggr) \frac{y_{\ast}}{y}\bigl(1+\alpha( x)
\bigr) \,dx.
\end{eqnarray*}
Recalling that $\alpha( x) =y^{2}/[x(x-y)y_{\ast}^{2}]$,
we can
use the fact that $y_{\ast}\geq\delta_{0}b$ to conclude that if
$x\in\lbrack b,b+a/b]$, then $\alpha( x) \leq\delta
_{0}^{-2}$, and
therefore
\[
\int_{b}^{b+a/b}\biggl( \frac{xy_{\ast}}{y}\bigl(1+2\alpha(
x)
\bigr)+1\biggr) \frac{y_{\ast}}{y}\bigl(1+\alpha( x) \bigr)\, dx\leq
4\delta
_{0}^{-8}a.
\]
We thus obtain that%
%
%e6.17 ###
%
\begin{equation}\label{IneqAux2}%
\frac{\int_{b}^{b+a/b}F^{\prime}(x) \,dx}{P\{ \sup_{t\in T}%
f(t)>b\} }\leq4a\delta_{0}^{-8}\frac{P\{ \sup_{t\in
T}f(t)>b-1/b\} }{P\{ \sup_{t\in T}f(t)>b\} }
\end{equation}
for any $b\geq b_{0}$. This inequality, together with the fact that
$F $ is continuously differentiable on $(-\infty,\infty)$,
yields the proof of the lemma for $b\geq0$.
\end{pf}

The previous result translates a question that involves the conditional
distribution of
$\max_{t\in T}f( t) $ near $b$ into a question involving the
tail distribution of $\max_{t\in T}f( t) $. The next
result then
provides a bound on this tail distribution.
\begin{lemma}
\label{LemBndDens1}For each $v>0$ there exists a constant $C(v)\in
(0,\infty)$
(possibly depending on $v>0$ but otherwise independent of $b$) so that such
that%
\[
P\Bigl( \max_{t\in T}f( t) >b\Bigr) \leq
C(v)b^{2d/\beta
+dv+1}\max_{t\in T}P\bigl( f( t) >b\bigr)
\]
for all $b\geq1$.
\end{lemma}
\begin{pf}
The proof of this result follows along the same lines of Theorem~2.6.2
in~\cite{ADT10}. Consider an open cover of $T=\bigcup_{i=1}^{M}T_{i}(\theta
)$, where
$T_{i}(\theta)=\{s\dvtx|s-t_{i}|<\theta\}$. We choose $t_{i}$
carefully such that
$N(\theta)=O(\theta^{-d})$ for $\theta$ arbitrarily small. Write
$f(
t) =g( t) +\mu( t) $, where $g(
t) $ is a centered Gaussian random field and note, using (A2) and
(A3), that
\[
P\Bigl( \max_{t\in T_{i}( \theta) }f( t)
>b\Bigr)
\leq P\Bigl( \max_{t\in T_{i}( \theta) }g(
t)
>b-\mu( t_{i}) -\kappa_{H}\theta^{\beta}\Bigr) .
\]
Now we\vspace*{1pt} wish to apply the Borel--TIS inequality (Theorem~\ref{ThmA2})
with $\mathcal{U}=T_{i}( \theta) $, $f_0=g$, $d( s,t) =E^{1/2}([ g( t)
-g( s) ] ^{2})$, which, as a consequence of (A2) and (A3), is bounded
above by $C_{0}\vert t-s\vert^{\beta/2}$ for some $C_{0}\in(0,\infty)$.
Thus, applying Theorem~\ref{ThmA1}, we have that $E\max_{t\in T_{i}(
\theta) }g( t) \leq
C_{1}%
\theta^{\beta/2}\log(1/\theta)$ for some $C_{1}\in(0,\infty)$.
Consequently,
the Borel--TIS inequality yields that there exists $C_{2}( v)
\in(0,\infty)$ such that for all $b$ sufficiently large and $\theta$
sufficiently small we have
\[
P\Bigl( \max_{t\in T_{i}( \theta) }g( t)
>b-\mu( t_{i}) -\kappa_{H}\theta^{\beta}\Bigr) \leq
C_{2}( v) \exp\biggl( -\frac{( b-\mu(
t_{i})
-C_{1}\theta^{\beta/(2+\beta v)}) ^{2}}{2\sigma
_{T_{i}}^{2}}\biggr) ,
\]
where $\sigma_{T_{i}}=\max_{t\in T_{i}(\theta)}\sigma(
t) $. Now
select $v>0$ small enough, and set $\theta^{\beta/(2+\beta v)}=b^{-1}$.
% which
%s $\theta=b^{-2/\beta-v}$ and conclude after
Straightforward calculations yield that%
\begin{eqnarray*}
P\Bigl( \max_{t\in T_{i}( \theta) }f( t)
>b\Bigr)
&\leq& P\Bigl( \max_{t\in T_{i}( \theta) }g(
t)
>b-\mu( t_{i}) -\kappa_{H}\theta^{\beta}\Bigr)
\\ &\leq&
C_{3}( v) \max_{t\in T_{i}( \theta) }\exp
\biggl(
-\frac{( b-\mu( t) ) ^{2}}{2\sigma(
t)
^{2}}\biggr)
\end{eqnarray*}
for some $C_{3}(v)\in(0,\infty)$. Now, recall the well-known
inequality (valid
for $x>0$) that
\[
\phi( x) \biggl( \frac{1}{x}-\frac{1}{x^{3}}\biggr)
\leq
1-\Phi( x) \leq\frac{\phi( x) }{x},\vadjust{\goodbreak}
\]
where $\phi=\Phi' $ is the standard Gaussian density.
Using this inequality it follows
that $C_{4}( v) \in( 0,\infty) $ can be
chosen so
that%
\[
\max_{t\in T_{i}( \theta) }\exp\biggl( -\frac{(
b-\mu( t) ) ^{2}}{2\sigma( t)
^{2}}\biggr)
\leq C_{4}( v) b\max_{t\in T_{i}}P\bigl( f(
t)
>b\bigr)
\]
for all $b\geq1$. We then conclude that there exists $C( v)
\in( 0,\infty) $ such that%
\begin{eqnarray*}
P\Bigl( \max_{t\in T}f( t) >b\Bigr) & \leq & N(
\theta) C_{4}( v) b\max_{t\in T}P\bigl( f(
t) >b\bigr) \\
& \leq & C\theta^{-d}b\max_{t\in T}P\bigl( f( t) >b\bigr)\\
&=&Cb^{2d/\beta+dv+1}\max_{t\in T}P\bigl( f( t) >b
\bigr)
\end{eqnarray*}
giving the result.
\end{pf}

We can now complete the proof of Proposition~\ref{PropDensity}.
\begin{pf*}{Proof of Proposition~\ref{PropDensity}}
The result is a straightforward
corollary of the previous two lemmas. By~(\ref{LemPart2}) in Lemmas
\ref{LemDens1} and~\ref{LemBndDens1} there exists
$\lambda\in(0,\infty)$ for which%
\begin{eqnarray*}
&&P\Bigl( \max_{t\in T}f(t)\leq b+a/b \big| \max_{t\in T}f(t)>b\Bigr)
\\
&&\qquad \leq
aA\frac{P( \max_{t\in T}f( t) \geq b-1/b)
}{P(
\max_{t\in T}f( t) \geq b) }\\
&&\qquad \leq aCAb^{2d/\beta+dv+1}\frac{\max_{t\in T}P
( f( t)
>b-1/b) }{P( \max_{t\in T}f( t) \geq b
) }\\
&&\qquad \leq aCAb^{2d/\beta+dv+1}\frac{\max_{t\in T}P
( f( t)
>b-1/b) }{\max_{t\in T}P( f( t) >b)
}\\
&&\qquad \leq a\lambda b^{2d/\beta+dv+1}.
\end{eqnarray*}
The last two inequalities follow from the obvious bound
\[
P\Bigl( \max_{t\in
T}f( t) \geq b\Bigr) \geq\max_{t\in T}P\bigl( f(
t) >b\bigr)
\]
and standard properties of the Gaussian
distribution. This yields~(\ref{PropPart1}), from which the remainder
of the proposition
follows.
\end{pf*}

%s7 ###
\section{Fine tuning: Twice differentiable homogeneous fields}
\label{SecSmooth}

In the preceding section we constructed a polynomial time algorithm
based on a
randomized discretization scheme. Our goal in this section is to illustrate
how to take advantage of additional information
to further improve the running time and the efficiency of the algorithm.
In order to illustrate our techniques we shall perform a more refined analysis
in the setting of smooth and homogeneous fields and shall establish
optimality of the algorithm in a~precise sense, to described below. Our
assumptions throughout this section are (B1) and (B2) of Section~\ref{SecMain}.

Let $C(s-t) = \Cov(f(s),f(t))$ be the covariance function of $f$,
which we assume also has mean zero. Note that it is an immediate
consequence of homogeneity and
differentiability that $\partial_{i} C(0)=\partial^{3}_{ijk}
C(0)=0$.\vspace*{1pt}

We shall need the following definition.
\begin{definition}
We call $\widetilde{T}=\{t_{1},\ldots,t_{M}\}\subset T$ a $\theta$-regular
discretization of $T$ if, and only if,
\[
\min_{i\neq j}|t_{i}-t_{j}|\geq\theta,\qquad \sup_{t\in T}\min
_{i}|t_{i}%
-t|\leq2\theta.
\]
\end{definition}

Regularity ensures that points in the grid $\widetilde{T}$ are well separated.
Intuitively, since $f$ is smooth, having tight clusters of points
translates to a
waste of computing resources, as a result of sampling highly correlated
values of $f$.
Also, note that every region containing a ball of radius $2\theta$ has at
least one representative in $\widetilde{T}$. Therefore, $\widetilde
{T}$ covers the
domain $T$ in an economical way. One technical convenience of $\theta
$-regularity is that for subsets $A\subseteq T$ that have positive Lebesgue
measure (in particular ellipsoids)
\[
\lim_{M\rightarrow\infty}\frac{\#(A\cap\widetilde{T})}{M}=\frac
{m(A)}{m(T)},
\]
where here and throughout the remainder of the section $\#(A)$ denotes the
cardinality of the set $A$.

Let $\widetilde{T}=\{t_{1},\ldots,t_{M}\}$ be a $\theta$-regular discretization
of $T$, and
consider
\[
{X}=( X_{1},\ldots,X_{M}) ^{T}\definedas( f(
t_{1}) ,\ldots,f( t_{M}) )^T .
\]
We shall concentrate on estimating $w_{M}( b) =P(
\max_{1\leq i\leq M}X_{i}>b) $. The next result (which we prove
in Section~\ref{SubBiasSmooth}) shows that if $\theta=\varepsilon
/b$, then the
relative bias is $O( \sqrt{\varepsilon}) $.
\begin{proposition}
\label{PropRelBias} Suppose $f$ is a Gaussian random field satisfying
conditions \textup{(B1)} and \textup{(B2)}. There exist $c_{0}$, $c_{1}$, $b_{0}$ and
$\varepsilon_{0}$ such that, for any finite $\varepsilon/b$-regular
discretization $\widetilde{T}$ of $T$,
%
%e7.1 ###
%
\begin{equation}\label{EqPrRelB}\quad
P\Bigl(\sup_{t\in\widetilde{T}}f(t)<b\big|\sup_{t\in T}f(t)>b\Bigr)\leq
c_{0}\sqrt{\varepsilon
}\quad \mbox{and}\quad \#(\widetilde{T})\leq c_{1}m(T)\varepsilon^{-d}b^{d}
\end{equation}
for all $\varepsilon\in(0,\varepsilon_{0}]$ and $b>b_{0}$.
\end{proposition}

Note that the bound on the bias obtained for twice differentiable
fields is
much sharper than that of the general H\"{o}lder\vadjust{\goodbreak} continuous fields
given by~(\ref{bias}) in Theorem~\ref{PropFinalHolder}.
This is partly because the conditional distribution of the random
field around local maxima is harder to describe in the H\"{o}lder continuous
than in the case of twice differentiable fields. In addition to the
sharper description
of the bias, we shall also soon show in Theorem~\ref{PropOpt} that our choice
of discretization is optimal in a cetain sense. Finally, we point out
that the
bound of $\sqrt{\varepsilon}$ in the first term of~(\ref{EqPrRelB})
is not
optimal. In fact, there seems to be some room of improvement, and we
believe that
a more careful analysis might yield a bound of the form
$c_{0}\varepsilon^{2}$.

We shall estimate $w_{M}( b) $ by using a slight variation
of Algorithm~\ref{algorithm1}. In particular, since the $X_{i}$'s are
now identically
distributed, we redefine $Q $ to~be
%
%e7.2 ###
%
\begin{equation} \label{SampHomo}%
Q( {X}\in B) =\sum_{i=1}^{M}\frac{1}{M}P[
{X}\in B\vert X_{i}>b-1/b].
\end{equation}
Our estimator then takes the form
%
%e7.3 ###
%
\begin{equation}\label{ImpEst}%
\widetilde{L}_{b}=\frac{M\times P(X_{1}>b-1/b)}{\sum_{j=1}^{M}\indic
(X_{j}>b-1/b) }\indic\Bigl(\max_{1\leq i\leq M}X_{i}>b\Bigr).
\end{equation}
Clearly, we have that $E^{Q}(L_{b})=w_{M}( b) $. (The reason
for subtracting the factor of $1/b$ was explained in Section~\ref{SubCE}.)
\begin{algorithm}
\label{algorithm3}
Given a number of replications $n$
and an $\varepsilon/b$-regular discretization $\widetilde{T}$
the algorithm is as follows:

\textsc{Step} (1). Sample ${X}^{( 1)
},\ldots,{X}^{(
n) }$ i.i.d. copies of ${X}$ with distribution $Q$ given by
(\ref{SampHomo}).

\textsc{Step} (2). Compute and output
\[
\widehat{L}_{n}=\frac{1}{n}\sum_{i=1}^{n}\widetilde{L}_{b}^{(
i) },
\]
where
\[
\widetilde{L}_{b}^{( i) }=\frac{M\times
P(X_{1}>b-1/b)}{\sum
_{j=1}^{M} \indic(X_{j}^{(i)}>b-1/b)} \indic\Bigl(\max_{1\leq i\leq
M}X_{j}^{( i)
}>b\Bigr).
\]
\end{algorithm}

Theorem~\ref{ThmBdRelHom} later guides the selection of $n$ in order to
achieve a prescribed relative error. In particular, our analysis, together
with\vspace*{1pt} considerations from Section~\ref{SecSet}, implies that choosing $n=O(
\varepsilon^{-2}\delta^{-1}) $ suffices to achieve
$\varepsilon$
relative error with probability at least $1-\delta$.

Algorithm~\ref{algorithm3} improves on Algorithm~\ref{algorithm2} for
H\"older continuous fields in two important ways.
The first aspect is that it is possible to obtain information on the
size of the relative bias of the estimator. In Proposition~\ref{PropRelBias},
we saw that in order to overcome bias due to discretization, it
suffices to
take a discretization of size $M=\#(\widetilde{T})=\Theta(b^{d})$.
That this
selection is also asymptotically optimal, in the sense
described in the next result, will be proven in Section~\ref{SubBiasSmooth}.\vadjust{\goodbreak}
\begin{theorem}
\label{PropOpt} Suppose $f$ is a Gaussian random field satisfying
conditions \textup{(B1)} and \textup{(B2)}. If $\theta\in(0,1)$, then, as $b\to\infty$,
\[
\sup_{\#(\widetilde{T})\leq b^{\theta d}}P\Bigl(\sup_{t\in\widetilde
{T}}f(t)>b\big|\sup_{t\in
T}f(t)>b\Bigr)\rightarrow0.
\]
\end{theorem}

This result implies that the relative bias goes to 100\% as
$b\to\infty$ if one chooses a discretization scheme of size $O(
b^{\theta d}) $ with $\theta\in(0,1)$. Consequently,
$d$ is the
smallest power of $b$ that achieves any given bounded relative bias,
and so
the suggestion above of choosing $M=O(b^d) $ points for
the discretization is,
in this sense, optimal.

The second aspect of improvement involves the variance. In the case of
H\"older continuous fields, the ratio of the second moment of the estimator
and $w( b) ^{2}$ was shown to be bounded by a quantity
that is
of order $O(M^{2})$. In contrast, in the context of smooth and homogeneous
fields considered here, the next result shows that this ratio is bounded
uniformly for $b>b_{0}$ and $M=\#(\widetilde{T})\geq cb^{d}$. That is, the
variance remains strongly controlled.
\begin{theorem}
\label{ThmBdRelHom}Suppose $f$ is a Gaussian random field satisfying
conditions \textup{(B1)} and \textup{(B2)}. Then there exist constants $c$, $b_{0}$ and
$\varepsilon_{0}$
such that for any $\varepsilon/b$-regular discretization $\widetilde
{T}$ of $T$
we have
\[
\sup_{b>b_{0},\varepsilon\in\lbrack0,\varepsilon_{0}]}\frac
{E^{Q}\widetilde{L}%
_{b}^{2}}{P^{2}( \sup_{t\in T}f( t) >b)
}\leq
\sup_{b>b_{0},\varepsilon\in\lbrack0,\varepsilon_{0}]}\frac
{E^{Q}\widetilde{L}
_{b}^{2}}{P^2( \sup_{t\in\widetilde{T}}f( t)
>b) }\leq c
\]
for some $c\in(0,\infty)$.
\end{theorem}

The proof of this result is given in Section~\ref{SubVarSmooth}. The
fact that
the number of replications remains bounded in $b$ is a consequence of the
strong control on the variance.

Finally, we note that the proof of Theorem~\ref{ThmHom} follows as a direct
corollary of Theorem~\ref{ThmBdRelHom} together with Proposition
\ref{PropRelBias} and our discussion in Section~\ref{SecSet}.
Assuming that
placing each point in $\widetilde{T}$ takes no more than $\mathbf{c}$
units of computer
time, the total complexity\vspace*{1pt} is, according to the discussion in Section
\ref{SecSet}, $O( nM^{3}+M) =O( \varepsilon
^{-2}\delta
^{-1}M^{3}+M) $. The contribution of the term $M^{3}=O(
\varepsilon^{-6d}b^{3d}) $ comes from the complexity of
applying Cholesky
factorization, and the term $M=O(\varepsilon^{-2d}b^d)$ corresponds to
the complexity
of placing $\widetilde{T}$.
\begin{remark}
Condition (B2) imposes a convexity assumption on the boundary of
$T$. This assumption, although convenient in the development of
the proofs of Theorems~\ref{PropOpt} and~\ref{ThmBdRelHom}, is
not necessary. The results can be generalized, at the expense
of increasing the length and the burden in the technical
development, to the case in which $T$ is a $d$-dimensional
manifold satisfying the so-called Whitney conditions
\cite{AdlerTaylor07}.
\end{remark}

The remainder of this section is devoted to the proof of Proposition
\ref{PropRelBias}, Theorems~\ref{PropOpt} and~\ref{ThmBdRelHom}.
%s7.1 ###
\subsection{\texorpdfstring{Bias control: Proofs of Proposition \protect\ref{PropRelBias} and Theorem \protect\ref{PropOpt}}
{Bias control: Proofs of Proposition 7.2 and Theorem 7.4}}\label{SubBiasSmooth}

We start with some useful lemmas, for all of which we assume that
Conditions (B1) and (B2) are satisfied. We shall also assume that the
global maximum of $f$ over $T$ is achieved, with probability one,
at a single point in $T$. Additional conditions under which this will happen
can be found in~\cite{AdlerTaylor07} and require little more than the
nondegeneracy of the joint distribution of $f$ and its first- and
second-order
derivatives. Of these lemmas, Lemma~\ref{LemLocalMaxima}, the proof of
which we defer to Section~\ref{SecProof}, is central to much of what follows.
However, before we state it we take a moment to describe Palm measures,
which may
not be familiar to all readers.

%s7.1.1 ###
\subsubsection{Palm distributions and conditioning}
\label{extentpalmsec}

It is well known that one needs to be careful treating the
distributions of
stochastic processes at random times. For a simple example, in the
current setting,
consider the behavior of a smooth stationary Gaussian process $f$ on
$\real$ along with
its derivative~$f'$. If $t\in\real$ is a fixed point, $u>0$, and we
are given that $f(0)=0$ and $f(t)=u$, then the conditional distribution
of~$f'(t)$ is
still Gaussian, with parameters determined by the trivariate
distribution of
$(f(0),f(t),f'(t))$. However, if we are given that $f(0)=0$, and that
$t>0$ is the
\textit{first} positive time that $f(t)=u$, then $t$ is an upcrossing of
the level
$u$ by $f$, and so $f'(t)$ must be positive. Thus it cannot be
Gaussian. The
difference between the two cases lies in the fact that in the first
case $t$ is
deterministic, while in the second it is random.

We shall require something similar, conditioning on the behavior of our
(Gaussian)
random
fields in the neighborhood of local maxima. Since local maxima are random
points, given their positions the distribution of the field is no longer
stationary nor, once again, even Gaussian. We often shall assume for
that a
local maximum is at the origin. This, however, amounts to saying that
the point-process induced by the set of local
maxima is Palm stationary (as opposed to space stationary) and
therefore we must then use the associated Palm distribution; the
precise conditional distribution of the field given the value of the
local maximum at the origin is given in Lemma~\ref{LemDist}.
The precise distribution is given in Lemma~\ref{LemDist}.

The theory behind this goes by the name of \textit{horizontal--vertical window
conditioning} and the resulting conditional distributions are known as
\textit{Palm distributions}. Standard treatments are given, for
example, in
\cite{GRF,ADT10,Kallenberg-Measures,Kallenberg-textbook,LLR}. To differentiate
between regular and Palm conditioning, we shall denote the latter by $\|
_{\mathcal P}$.

We can now set up two important lemmas which tell us
about the behavior of $f$ in the neighborhood of local and global
maxima. Proofs are deferred until Section~\ref{SecProof}. First, we
provide some
notation.

Let $\mathcal L$ be the (random) set of local maxima of $f$. That is,
for each
$s$ in the interior of $T$, $s\in\mathcal L$ if and only if
%
%e7.4 ###
%
\begin{equation}\label{cond}
\nabla f(s) = 0\quad \mbox{and}\quad \nabla^2 f(s) \in\mathcal N,
\end{equation}
where $\mathcal N$ is the set of negative definite matrices, and
$\nabla^2 f(s)$ is the Hessian matrix of $f$ at $s$. For $s\in
\partial T$, similar constraints apply and are
described in the proof of Lemma~\ref{LemLocalMaxima}.
Then we have:
\begin{lemma}
\label{LemLocalMaxima}
Let $\mathcal{L}$ be the set of local maxima of $f$.
For any $a_{0}>0$, there exists $c^{\ast}$,
$\delta^{\ast}$, $b_{0}$ and $\delta_{0}$ (which depend on the
choice of
$a_{0}$), such that for any $s\in\mathcal{L}$,
$a\in(0,a_{0})$, $\delta\in(0,\delta
_{0})$, $b>b_{0}$, $z> b+a/b$
%
%e7.5 ###
%
\begin{equation}
\label{Palm}
P\Bigl(\min_{|t-s|<\delta ab^{-1}}f(t)<b \big\|_{\mathcal P} f(s)=z\Bigr)
\leq
c^{\ast
}\exp\biggl( -\frac{\delta^{\ast}}{\delta^{2}}\biggr) .
\end{equation}
\end{lemma}
\begin{lemma}
\label{LemGlobalMaxima} Let $t^{\ast}$ be the point in $T$
at which the global maximum
of $f$ is attained.
Then, with the same choice of constants as in
Lemma~\ref{LemLocalMaxima}, for any $a\in(0,a_{0})$, $\delta\in
(0,\delta_{0}
)$ and $b>b_{0}$,
\[
P\Bigl(\min_{|t-t^{\ast}|<\delta ab^{-1}}f(t)<b
\big\|_{\mathcal P} f(t^{\ast})>b+a/b\Bigr)\leq2c^{\ast}
\exp\biggl( -\frac{\delta^{\ast}}{\delta^{2}}\biggr) .
\]
\end{lemma}

% The next lemma is a
%direct corollary of Lemma~\ref{LemLocalMaxima}.%
%
%at which the global maximum
%of $f$ is attained.
%Then, with the same choice of constants as in
%Lemma~\ref{LemLocalMaxima}, for any $a\in(0,a_{0})$, $\delta\in(0,
%)$, and $b>b_{0}$, we have
%P(\min_{|t-t^{\ast}|<\delta ab^{-1}}f(t)<b|f(t^{\ast})=b+a/b)\leq2c^{
%
%
%Note that, for any $s\in T$,
%& P( \min_{|t|\leq\delta a/b}f(t)<b|f(s)=b+a/b,s=t^{\ast})
%& =\frac{P( \min_{|t|\leq\delta a/b}f(t)<b,f(s)=b+a/b,s=t^{
%}) }{P( f(s)=b+a/b,s=t^{\ast}) }\\
%& \leq\frac{P( \min_{|t|\leq\delta a/b}f(t)<b,f(s)=b+a/b,s
%5& =P(\min_{|t-s|<\delta ab^{-1}}f(t)<b|f(s)=b+a/b,s\in
%)\frac{P( f(s)=b+a/b,s\in\mathcal{L}) }{P(
%f(s)=b+a/b,s=t^{\ast}) }.\\
%& \leq c^{\ast}\exp( -\frac{\delta^{\ast}}{\delta^{2}}
%)
%f(s)=b+a/b,s=t^{\ast}) }.%
%Applying this bound and the fact that
%f(s)=b+a/b,s=t^{\ast}) }\rightarrow1
%as $b\rightarrow\infty$, the conclusion of the Lemma holds.

%s7.1.2 ###
\subsubsection{Back to the proofs}
The following lemma gives a bound on the density of $\sup_{t\in
T}f(t)$, which
will be used to control the size of overshoot beyond level $b$.
\begin{lemma}
\label{LemOvershoot}Let $p_{f^{\ast}}(x)$ be the density function of
$\sup_{t\in T}f(t)$. Then there exists a constant $c_{f^{\ast}}$ and $b_{0}$
such that
\[
p_{f^{\ast}}(x)\leq c_{f^{\ast}}x^{d+1}P\bigl(f(0)>x\bigr)
\]
for all $x>b_{0}$.
\end{lemma}
\begin{pf}
Recalling (\ref
{mainresultintro}), let the continuous function
$p^{E}(x)$, $x\in\real$, be defined by the relationship
\[
E\bigl(\chi\bigl(\{t\in T\dvtx f(t)\geq b\} \bigr) \bigr) = \int_{b}^{\infty}p^{E}(x) \,dx,
\]
where the left-hand side is the expected value of the Euler--Poincar\'
{e} characteristic of $A_{b}$. Then,
according to Theorem 8.10 in~\cite{AzWschebor09}, there exists $c$ and~$\delta$ such that
\[
|p^{E}(x)-p_{f^{\ast}}(x)|<cP\bigl(f(0)>(1+\delta)x\bigr)
\]
for all $x>0$. In addition, thanks to the result of \cite
{AdlerTaylor07} which
provides $\int_{b}^{\infty}p^{E}(x)\,dx$ in closed form, there exists $c_{0}$
such that, for all $x>1$,
\[
p^{E}(x)<c_{0}x^{d+1}P\bigl(f(0)>x\bigr).\vadjust{\goodbreak}
\]
Hence, there exists $c_{f^{\ast}}$ such that
\[
p_{f^{\ast}}(x)\leq c_{0}x^{d+1}P\bigl(f(0)>x\bigr)+cP\bigl(f(0)>(1+\delta)x\bigr)\leq
c_{f^{\ast
}}x^{d+1}P\bigl(f(0)>x\bigr)
\]
for all $x>1$.
\end{pf}

The last ingredients required to provide the proof of Proposition
\ref{PropRelBias} and Theorem~\ref{PropOpt} are stated in the following
result, adapted from
Lemma 6.1 and Theorem~7.2 in~\cite{PiterbargGaussian} to the twice
differentiable case.
\begin{theorem}
There exists a constant $H$ (depending on the covariance
function $C$), such that
%
%e7.6 ###
%
\begin{equation}
\label{ThmLDHomo}
P\Bigl(\sup_{t\in T}f(t)>b\Bigr)=\bigl(1+o( 1) \bigr)Hm(T)b^{d}P\bigl(f(0)>b\bigr)
\end{equation}
as $b\rightarrow\infty$.

Similarly, %\label{LemPick}
choose $\delta$ small enough so that $[0,\delta]^{d}\subset
T$, and let $\Delta_{0} = [0,b^{-1}]^{d}$. Then there exists a constant
$H_{1}$ such that
%
%e7.7 ###
%
\begin{equation}
\label{LemPick}
P\Bigl(\sup_{t\in\Delta_{0}} f(t)> b\Bigr) = \bigl(1+o(1)\bigr) H_{1} P\bigl(f(0)>b\bigr)
\end{equation}
as $b\rightarrow\infty$.
\end{theorem}

%T$, and let $\Delta_{0} = [0,b^{-1}]^{d}$. Then there exists a constant
%$H_{1}$ such that
%P(\sup_{t\in\Delta_{0}} f(t)> b) = (1+o(1)) H_{1} P(f(0)>b),
%as $b\rightarrow\infty$.

We now are ready to provide the proof of Proposition~\ref{PropRelBias} and
Theorem~\ref{PropOpt}.
\begin{pf*}{Proof of Proposition~\ref{PropRelBias}}
The fact that there
exists $c_{1}$
such that
\[
\#(\widetilde{T})\leq c_{1}m(T)\varepsilon^{-d}b^{d}%
\]
is immediate from assumption (B2). Therefore, we proceed to provide a~bound for
the relative bias. Note first that elementary conditional probability
manipulations yield that, for any
$\varepsilon>0$,
\begin{eqnarray*}
&& P\Bigl(\sup_{t\in\widetilde{T}}f(t)<b\big|\sup_{t\in T}f(t)>b\Bigr)\\
&&\qquad \leq P\Bigl(\sup_{t\in T}f(t)<b+2\sqrt{\varepsilon}/b\big|\sup_{t\in T}%
f(t)>b\Bigr)\\
&&\qquad\quad{}+P\Bigl(\sup_{t\in\widetilde{T}}f(t)<b\big|\sup_{t\in T}f(t)>b+2\sqrt
{\varepsilon
}/b\Bigr).
\end{eqnarray*}
By~(\ref{ThmLDHomo}) and Lemma~\ref{LemOvershoot}, there exists
$c_{2}$ such that, for large enough $b$, the first term above can
bounded by
\[
P\Bigl(\sup_{t\in T}f(t)<b+2\sqrt{\varepsilon}/b\big|\sup_{t\in
T}f(t)>b\Bigr)\leq
c_{2}\sqrt{\varepsilon}.
\]
Now take $\varepsilon<\varepsilon_{0}<\delta_{0}^{2}$ where $\delta
_{0} $
is as in Lemmas~\ref{LemLocalMaxima} and~\ref{LemGlobalMaxima}.
Then, applying~(\ref{Palm}), the second term can be bounded by
\begin{eqnarray*}
&& P\Bigl(\sup_{t\in\widetilde{T}}f(t)<b\big|\sup_{t\in T}f(t)>b+2\sqrt
{\varepsilon}/b\Bigr)\\[-2pt]
&&\qquad \leq P\Bigl(\sup_{|t-t^{\ast}|<2\varepsilon
b^{-1}}f(t)<b\big|\sup_{t\in
T}f(t)>b+2\sqrt{\varepsilon}\Bigr)\\[-2pt]
&&\qquad \leq2c^{\ast}\exp( -\delta^{\ast
}\varepsilon^{-1}) .
\end{eqnarray*}
Hence, there exists a $c_{0}$ such that
\[
P\Bigl(\sup_{t\in\widetilde{T}}f(t)<b\big|\sup_{t\in T}f(t)>b\Bigr)\leq
c_{2}\sqrt{\varepsilon
}+2c^{\ast}e^{-\delta^{\ast}/\varepsilon}\leq c_{0}\sqrt
{\varepsilon}
\]
for all $\varepsilon\in(0,\varepsilon_{0})$.\vspace*{-2pt}
\end{pf*}
\begin{pf*}{Proof of Theorem~\ref{PropOpt}}
We write $\theta=1-3\delta
\in(0,1)$. First
note that, by~(\ref{ThmLDHomo}),
\[
P\Bigl(\sup_{T}f(t)>b+b^{2\delta-1}\big|\sup_{T}f(t)>b\Bigr)\rightarrow0
\]
as $b\to\infty$. Let $t^{\ast}$ be the position of the global maximum
of $f$ in $T$. According to the exact Slepian model
in Section~\ref{SecProof} and an argument similar to the proof of
Lemmas~\ref{LemLocalMaxima} and~\ref{LemGlobalMaxima}
%
%e7.8 ###
%
\begin{equation} \label{far}%
P\Bigl(\sup_{|t-t^{\ast}|>b^{2\delta-1}}f(t)>b\big|b<f(t^{\ast})\leq
b+b^{2\delta
-1}\Bigr)\rightarrow0
\end{equation}
as $b\to\infty$. Consequently,
\[
P\Bigl(\sup_{|t-t^{\ast}|>b^{2\delta-1}}f(t)>b\big|\sup
_{T}f(t)>b\Bigr)\rightarrow0.
\]
Let%
\[
B(\widetilde{T},b^{2\delta-1})=\bigcup_{t\in\widetilde
{T}}B(t,b^{2\delta-1}).
\]
We have
\begin{eqnarray*}
&& P\Bigl(\sup_{\widetilde{T}}f(t)>b\big|\sup_{T}f(t)>b\Bigr)\\[-2pt]
&&\qquad \leq P\Bigl(\sup_{|t-t^{\ast}|>b^{2\delta-1}}f(t)>b\big|\sup
_{T}f(t)>b\Bigr)\\[-2pt]
&&\qquad\quad{} +P\Bigl(\sup_{\widetilde{T}}f(t)>b,\sup
_{|t-t^{\ast}|>b^{2\delta-1}}f(t)\leq
b\big|\sup_{T}f(t)>b\Bigr)\\[-2pt]
&&\qquad \leq o(1)+P\Bigl(t^{\ast}\in B(\widetilde{T},b^{2\delta
-1})\big|\sup_{T}f(t)>b\Bigr)\\[-2pt]
&&\qquad \leq o(1)+P\Bigl(\sup_{B(\widetilde{T},b^{2\delta
-1})}f(t)>b\big|\sup_{T}f(t)>b\Bigr).\vadjust{\goodbreak}
\end{eqnarray*}
Since $\#(\widetilde{T})\leq b^{(1-3\delta)d}$, we can find a finite set
$T^{\prime}=\{t^{\prime}_{1},\ldots,t^{\prime}_{l}\} \subset T$ and let
$\Delta_{k} = t^{\prime}_{k} + [0,b^{-1}]$ such that\vspace*{2pt}
$l=O(b^{(1-\delta) d})$
and $B(\widetilde T, b^{2\delta-1 })\subset\bigcup_{k=1}^{l} \Delta
_{k}$. The choice
of $l$ only depends on $\#(\widetilde T)$, not the particular
distribution of
$\widetilde T$. Therefore, applying~(\ref{LemPick}),
\[
\sup_{\#(\widetilde T \leq b^{\theta d})}P\Bigl(\sup_{B(\widetilde
{T},b^{2\delta-1}%
)}f(t)>b\Bigr)\leq O\bigl(b^{(1-\delta) d}\bigr) P\bigl(f(0)>b\bigr).
\]
This, together with~(\ref{ThmLDHomo}), yields
\[
\sup_{\#(\widetilde T \leq b^{\theta d})}P\Bigl(\sup_{B(\widetilde
{T},b^{2\delta-1}%
)}f(t)>b\big|\sup_{T}f(t)>b\Bigr)\leq O(b^{-\delta d }) = o(1)
\]
for $b\geq b_{0}$, which clearly implies the statement of the result.
\end{pf*}

%s7.2 ###
\subsection{\texorpdfstring{Variance control: Proof of Theorem \protect\ref{ThmBdRelHom}}{Variance control: Proof of Theorem 7.5}}
\label{SubVarSmooth}

We proceed directly to the proof of Theorem~\ref{ThmBdRelHom}.
\begin{pf*}{Proof of Theorem~\ref{ThmBdRelHom}}
Note that
\begin{eqnarray*}
&& \frac{E^{Q}\widetilde{L}_{b}^{2}}{P^{2}( \sup_{t\in T}f
( t)
>b) }
\\ &&\qquad =\frac{E(\widetilde{L}_{b})}{P^{2}( \sup_{t\in
T}f(
t) >b) }\\
&&\qquad =\frac{E( {M\times P(X_{1}>b-1/b)}/{\sum
_{j=1}^{n}\indic(
X_{j}>b-{1/b}) };\max_{j}X_{j}>b) }{P^{2}(
\sup_{t\in
T}f( t) >b)} \\
&&\qquad =\biggl({E\biggl( \frac{M\times P(X_{1}>b-1/b)}{\sum
_{j=1}^{n}\indic(
X_{j}>b-{1/b}) };\max_{j}X_{j}>b,\sup_{t\in T}f(
t)
>b\biggr) }\biggr)\\
&&\qquad\quad{}\times\Bigl({P^{2}\Bigl( \sup_{t\in T}f( t) >b\Bigr)}\Bigr)^{-1}
\\
&&\qquad =\biggl({E\biggl( \frac{MP(X_{1}>b-1/b)\indic(
\max_{j}X_{j}>b)
}{\sum_{j=1}^{n}\indic( X_{j}>b-1/b) }\Big\vert\sup
_{t\in
T}f( t) >b\biggr)}\biggr)\\
&&\qquad\quad{}\times\Bigl({P\Bigl( \sup_{t\in T}f(
t)
>b\Bigr)}\Bigr)^{-1} \\
&&\qquad =E\biggl( \frac{M\indic( \max
_{j}X_{j}>b) }{\sum_{j=1}%
^{n}\indic( X_{j}>b-1/b) }\Big\vert\sup_{t\in
T}f( t)
>b\biggr) \\
&&\qquad\quad{}\times\frac{P(X_{1}>b-1/b)}{P( \sup_{t\in T}f(
t)
>b) }.
\end{eqnarray*}
The remainder of the proof involves showing that the last conditional
expectation
here is of order $O(b^{d})$. This, together with
(\ref{ThmLDHomo}), will yield the result. Note that for any $A(
b,\varepsilon) $ such that $A( b,\varepsilon)
=\Theta(M)$
uniformly over $b$ and $\varepsilon$, we can write
%
%e7.9 ###
%
\begin{eqnarray} \label{Decomp}
&& E\biggl( \frac{M\times\indic( \max_{j}X_{j}>b
) }%
{\sum_{j=1}^{M}\indic( X_{j}>b-{1/b}) }
\Big\vert\sup_{t\in
T}f( t) >b\biggr)\nonumber\\
&&\qquad \leq E\biggl( \frac{M\times\indic( \sum
_{j=1}^{M}\indic(
X_{j}>b-1/b) \geq{M}/{A( b,\varepsilon)
})
}{\sum_{j=1}^{M}\indic( X_{j}>b-{1/b}) }
\Big\vert\sup_{t\in
T}f( t) >b\biggr) \nonumber\\[-8pt]\\[-8pt]
&&\qquad\quad{} +E\biggl( \frac{M\times\indic( 1\leq
\sum_{j=1}^{M}\indic(
X_{j}>b-1/b) <{M}/{A( b,\varepsilon) }
) }%
{\sum_{j=1}^{M}\indic( X_{j}>b-{1/b}) }
\Big\vert\nonumber\\
&&\hspace*{232.3pt}\sup_{t\in
T}f( t) >b\biggr) .\nonumber
\end{eqnarray}
We shall select $A(b,\varepsilon)$ appropriately in order to bound the
expectations above. By Lemma~\ref{LemGlobalMaxima}, %
for any $4\varepsilon\leq\delta\leq\delta_{0}$, there exist constants
$c^{\prime}$ and $c^{\prime\prime}\in( 0,\infty) $,
such that
%
%e7.10 ###
%
\begin{eqnarray}\label{BDArea}\quad
c^{\ast}\exp\biggl( -\frac{\delta^{\ast}}{\delta^{2}}\biggr)
&\geq&P\Bigl( \min_{\vert t-t^{\ast}\vert\leq
\delta
/b}f( t) <b-1/b\big\vert\sup_{t\in T}f( t)
>b\Bigr) \nonumber\\
& \geq& P\Biggl( \sum_{j=1}^{M}\indic(
X_{j}>b-1/b) \leq
c^{\prime}\delta^{d}/\varepsilon^{d}\big\vert\sup_{t\in
T}f(
t) >b\Biggr)\\
& \geq& P\biggl( \frac{M}{\sum_{j=1}^{M}\indic
(X_{j}>b-1/b)}\geq\frac
{b^{d}c^{\prime\prime}}{\delta^{d}}\Big\vert\sup_{t\in T}f
( t)
>b\biggr) .\nonumber
\end{eqnarray}
The first inequality is an application of Lemma~\ref{LemGlobalMaxima}.
%%{LemLocalMaxima}
The
second inequality is due to the fact that for any ball $B$ of radius
$4\varepsilon$ or larger, $\#(\widetilde{T}\cap B)\geq c^{\prime
d}\varepsilon
^{-d}$ for some \mbox{$c^{\prime}>0$}. Inequality~(\ref{BDArea}) implies
that for all
$x$ such that $b^{d}c^{\prime\prime}/\delta_{0}^{d}<x<b^{d}c^{\prime
\prime
}/[4^{d}\varepsilon^{d}]$, there exists $\delta^{\ast\ast}>0$ such that
\[
P\biggl( \frac{M}{\sum_{j=1}^{M}\indic( X_{j}>b-
{1/b})
b^{d}}\geq x\big\vert\sup_{t\in T}f( t) >b\biggr)
\leq
c^{\ast}\exp( -\delta^{\ast\ast}x^{2/d}) .
\]
Now let $A( b,\varepsilon) =b^{d}c^{\prime\prime}%
/(4^{d}\varepsilon^{d})$ and observe that by the second result in
(\ref{EqPrRelB}) we have
$A( b,\varepsilon) =\Theta( M) $ and,
moreover, that
there exists $c_{3}$ such that the first term on the right-hand side of
(\ref{Decomp}) is bounded by
%
%e7.11 ###
%
\begin{eqnarray} \label{Large}
&&E\biggl( \frac{M\times\indic( \sum_{j=1}^{M}\indic
( X_{j}%
>b-1/b) \geq{M}/{A( b,\varepsilon) })
}%
{\sum_{j=1}^{M}\indic( X_{j}>b-1/b) }\Big\vert\sup
_{T}f(
t) >b\biggr) \nonumber\\[-8pt]\\[-8pt]
&&\qquad\leq c_{3}b^{d}.\nonumber
\end{eqnarray}
Now we turn to the second term on the right-hand side of (\ref
{Decomp}). We use the fact that
$M/A(b,\varepsilon)\leq c^{\prime\prime\prime}\in(0,\infty)$
(uniformly as
$b\rightarrow\infty$ and $\varepsilon\rightarrow0$). There exist $c_{4}$
and $c_{5}$ such that, for $\varepsilon\leq\delta_{0}/c_{4}$,
%
%e7.12 ###
%
\begin{eqnarray} \label{Small}
&& E\biggl( \frac{M\indic( 1\leq\sum_{j=1}^{M}\indic
( X_{j}%
>b-{1}/{b}) <{M}/{A( b,\varepsilon)
})
}{\sum_{j=1}^{M}\indic( X_{j}>b-{1}/{b}) }
\Big\vert\sup
_{T}f( t) >b\biggr)\nonumber\\
&&\qquad \leq MP\Biggl( \sum_{j=1}^{n}\indic\biggl(
X_{j}>b-\frac{1}{b}\biggr)
<c^{\prime\prime\prime}\big\vert\sup_{T}f( t)
>b\Biggr)
\nonumber\\[-8pt]\\[-8pt]
&&\qquad \leq MP\Bigl( \min_{\vert t-t^{\ast
}\vert\leq
c_{4}\varepsilon/b}f( t) <b-1/b\big\vert\sup
_{T}f(
t) >b\Bigr) \nonumber\\
&&\qquad \leq c_{1}m(T)b^{d}\varepsilon^{-d}\exp\biggl(
-\frac{\delta^{\ast}}%
{c_{4}^{2}\varepsilon^{2}}\biggr) \leq c_{5}b^{d}.\nonumber
\end{eqnarray}
The second inequality holds from the fact that if $\sum
_{j=1}^{M}\indic( X_{j}%
>b-\frac{1}{b}) $ is less than~$c^{\prime\prime\prime}$,
then the minimum of
$f $ in a ball around the local maximum and of radius
$c_{4}\varepsilon/b$ must be less than $b-1/b$. Otherwise, there are
more than~$c^{\prime\prime\prime}$ elements of $\widetilde{T}$ inside such a
ball. The last
inequality is due to Lemma~\ref{LemGlobalMaxima} and Theorem \ref
{PropRelBias}.

Putting~(\ref{Large}) and~(\ref{Small}) together we obtain, for all
$\varepsilon/b$-regular discretizations with $\varepsilon
<\varepsilon_{0}%
=\min(1/4,1/c_{4})\delta_{0}$,
\begin{eqnarray*}
\frac{E^{Q}\widetilde{L}_{b}^{2}}{P( \sup_{t\in T}f(
t)
>b) } & \leq & E\biggl( \frac{M\indic( \max
_{j}X_{j}>b)
}{\sum_{j=1}^{M}\indic( X_{j}>b-{1/b}) }
\Big\vert\sup_{t\in
T}f( t) >b\biggr) \\
&&{}\times\frac{P(X_{1}>b-1/b)}{P( \sup
_{t\in
T}f( t) >b) }\\
& \leq &(c_{3}+c_{5})\frac{b^{d}P(X_{1}>b-1/b)}{P( \sup_{t\in
T}f(
t) >b) }.
\end{eqnarray*}
Applying now~(\ref{ThmLDHomo})
%P(\sup_{T}f(t)>b)=(1+o(1))Hm(T)b^{d}P(X_{1}>b),
and Proposition~\ref{PropRelBias}, we have that
\[
P\Bigl(\sup_{t\in T}f(t)>b\Bigr)<\frac{P(\sup_{t\in\widetilde
{T}}f(t)>b)}{1-c_{0}%
\sqrt{\varepsilon}},
\]
and we have
\[
\sup_{b>b_{0},\varepsilon\in\lbrack0,\varepsilon_{0}]}\frac
{E^{Q}\widetilde{L}%
_{b}^{2}}{P( \sup_{t\in\widetilde{T}}f( t)
>b) }<\infty
\]
as required.
\end{pf*}

%s7.3 ###
\subsection{Remaining proofs} \label{SecProof}
We start with the proof of
Lemma~\ref{LemLocalMaxima}.
Without loss of generality, we assume that the random field
of that result has mean zero and unit
variance. However,\vadjust{\goodbreak} before getting into the details of the proof of
Lemma~\ref{LemLocalMaxima}, we need a
few additional lemmas, for which we adopt the following notation:
Let $C_{i}$ and $C_{ij} $ be the first- and second-order
derivatives of $C$, and
define the vectors
\begin{eqnarray*}
\mu_{1}(t) & = &(-C_{1}(t),\ldots,-C_{d}(t)),\\
\mu_{2}(t) & = &\operatorname{vech}\bigl((C_{ij}(t),i=1,\ldots,d,j=i,\ldots,d)\bigr).
\end{eqnarray*}
Let $f^{\prime}(0)$ and $f^{\prime\prime}(0)$ be the gradient and
vector of
second-order derivatives of $f$ at~$0$, where $f^{\prime\prime}(0)$
is arranged in the same
order as $\mu_{2}(0)$. Furthermore, let $\mu_{02}=\mu_{20}^{\top}$
be a vector
of second-order spectral moments and $\mu_{22}$ a matrix of fourth-order
spectral moments. The vectors $\mu_{02}$ and $\mu_{22}$ are arranged so
that
\[
\pmatrix{
1 & 0 & \mu_{02}\cr
0 & \Lambda& 0\cr
\mu_{20} & 0 & \mu_{22}}
\]
is the covariance matrix of $(f(0),f^{\prime}(0),f^{\prime\prime}(0))$,
where $\Lambda=(-C_{ij}(0))$. It then follows that
\[
\mu_{2\cdot0}=\mu_{22}-\mu_{20}\mu_{02}%
\]
be the conditional variance of $f^{\prime\prime}(0)$ given $f(0)$. The
following lemma, given in~\cite{ADT10}, provides a stochastic representation
of the $f$ given that it has a local maxima at level $u$ at the origin.
We emphasize that, as described above, the conditioning here is in the sense
of Palm distributions. The resultant conditional, or ``model'' process
(\ref{ModelProcess}) is generally called a \textit{Slepian process}.
\begin{lemma}
\label{LemDist} Given that $f$ has a local maximum with height $u$ at
zero (an
interior point of $T$), the conditional field is equal in distribution to
%
%e7.13 ###
%
\begin{equation}
\label{ModelProcess} \label{Slepian}%
f_{u}(t)\triangleq uC( t) -W_{u}\beta^{\top}(
t)
+g( t) .
\end{equation}
$g(t)$ is a centered Gaussian random field with covariance function
\[
\gamma(s,t)=C(s-t)-(C(s),\mu_{2}(s))
\pmatrix{
1 & \mu_{02}\cr
\mu_{20} & \mu_{22}}^{-1}
\pmatrix{
C(t)\cr
\mu_{2}^{\top}(t)}
-\mu_{1}(s)\Lambda^{-1}\mu_{1}^{\top}(t),
\]
and $W_{u}$ is a $\frac{d(d+1)}{2}$ random vector independent of
$g(t)$ with
density function
%
%e7.14 ###
%
\begin{equation} \label{DenW}%
\psi_{u}( w) \propto\bigl|{\det}\bigl( r^{\ast}(w)-u\Lambda
\bigr)
\bigr|\exp\bigl( -\tfrac{1}{2}w^{\top}\mu_{2\cdot0}^{-1}w\bigr)
\indic\bigl( r^{\ast
}(w)-u\Lambda\in\mathcal{N}\bigr) ,\hspace*{-28pt}
\end{equation}
where $r^{\ast}(w)$ is a $d\times d$ symmetric matrix whose upper triangular
elements consist of the components of $w$. The set of negative definite
matrices is denoted by $\mathcal{N}$. Finally, $\beta(t)$ is defined by
\[
( \alpha( t) ,\beta( t) )
=(
C( t) ,\mu_{2}( t) )
\pmatrix{
1 & \mu_{02}\cr
\mu_{20} & \mu_{22}%
}^{-1}.
\]
\end{lemma}

The following two technical lemmas, which we shall prove after completing
the proof of Lemma~\ref{LemLocalMaxima},
provide bounds for the last two terms of
(\ref{Slepian}).\vadjust{\goodbreak}
\begin{lemma}
\label{LemBeta}Using the notation in~(\ref{Slepian}), there exist
$\delta
_{0}$, $\varepsilon_{1}$, $c_{1}$ and $b_{0}$ such that, for any $u>b>b_{0}$
and $\delta\in(0,\delta_{0})$,
\[
P\biggl( \sup_{|t|\leq\delta a/b}\vert W_{u}\beta^{\top}
( t)
\vert>\frac{a}{4b}\biggr) \leq c_{1}\exp\biggl( -\frac
{\varepsilon
_{1}b^{2}}{\delta^{4}} \biggr) .
\]
\end{lemma}
\begin{lemma}
\label{LemG}There exist $c$, $\widetilde{\delta}$ and $\delta_{0}$
such that, for
any $\delta\in(0,\delta_{0})$,
\[
P\biggl( \max_{|t|\leq\delta a/b}|g(t)|>\frac{a}{4b}\biggr) \leq
c\exp\biggl(
-\frac{\widetilde{\delta}}{\delta^{2}}\biggr) .
\]
\end{lemma}
\begin{pf*}{Proof of Lemma~\ref{LemLocalMaxima}}
Using the notation of Lemma
\ref{LemDist}, given any \mbox{$s\in\mathcal{L}$} for which $f(s)=b$, we
have that
the corresponding conditional distribution of $s$ is that of
%f_{b}(\cdot)\equalinlaw
$f_b(\cdot-s)$.
Consequently, it suffices to show that
\[
P\Bigl( \min_{\vert t\vert\leq\delta
a/b}f_{b}(t)<b-a/b\Bigr)
\leq c^{\ast}\exp\biggl( -\frac{\delta^{\ast}}{\delta^{2}}
\biggr) .
\]

We consider first the case for which the local maximum is
in the interior of $T$. Then, by the Slepian model~(\ref{Slepian}),
\[
f_{b}( t) =bC( t) -W_{b}\beta^{\top}(
t) +g( t) .
\]
We study the three terms of the Slepian model individually. Since
\[
C(t)=1-t^{\top}\Lambda t+o( \vert t\vert^{2}
) ,
\]
there exists a $\varepsilon_{0}$ such that
\[
bC( t) \geq b-\frac{a}{4b}
\]
for all $|t|<\varepsilon_{0}\sqrt{a}/b$. According to Lemmas \ref
{LemBeta} and
\ref{LemG}, for $\delta<\min(\varepsilon_{0}/\sqrt{a_{0}}$, $\delta
_{0}) $,
\begin{eqnarray*}
&& P\Bigl( \min_{|t|\leq\delta a/b}f_{u}(t)<b-a/b\Bigr) \\
&&\qquad \leq P\biggl( \max_{|t|<\delta a/b}\vert g
( t)
\vert>\frac{a}{4b}\biggr) +P\biggl( \sup_{\vert
t\vert
\leq\delta a/b}\vert W_{b}\beta^{\top}(t)\vert>\frac
{a}{4b}\biggr)
\\
&&\qquad \leq c\exp\biggl( \frac{\widetilde{\delta}}{\delta
^{2}}\biggr) +c_{1}
\exp\biggl( -\frac{\varepsilon_{1}b^{2}}{\delta^{4}}\biggr) \\
&&\qquad \leq c^{\ast}\exp\biggl( -\frac{\delta^{\ast
}}{\delta^{2}}\biggr)
\end{eqnarray*}
for some $c^{\ast}$ and $\delta^{\ast}$.

Now consider the case for which the local maximum is in the $(d-1)$-dimensional
boundary of $T$. Due to convexity of $T$ we can assume, without loss of
generality,
that the tangent space of $\partial T$ is generated by
$\partial/\partial t_{2},\ldots,\break\partial/\partial t_{d}$,
the local maximum\vadjust{\goodbreak} is located at the origin and $T$ is a subset of the
positive half-plane $t_1\geq0$. That these arguments to not involve a
loss of
generality follows from the arguments on pages 192--291 of
\cite{AdlerTaylor07}, which rely on the assumed stationarity of $f$
(for translations) and the fact that rotations, while changing the
distributions, will not
affect the probabilities that we are currently computing.

For the origin, positioned as just described, to be a local maximum it is
necessary and sufficient that the gradient of
$f$ restricted to $\partial T$ is the zero vector, the Hessian matrix restricted
to $\partial T$ is negative definite and $\partial_{1}f(0)\leq0$. Applying
a version of Lemma~\ref{LemDist} for this case, conditional on
$f(0)=u$ and
$0$ being a local maximum, the field is equal in distribution to
%
%e7.15 ###
%
\begin{equation}\label{SlepianB}%
uC(t)-\widetilde{W}_{u}\beta^{\top}(t)+\mu_{1}(t)\Lambda^{-1}(
Z,0,\ldots,0) ^{T}+g(t),
\end{equation}
where $Z\leq0$ corresponds to $\partial_{1}f( 0) $ and
it follows
a truncated (conditional on the negative axis) Gaussian random variable with
mean zero and a variance parameter which is computed as the conditional
variance of $\partial_{1}f( 0) $ given $(\partial
_{2}f(
0) ,\ldots,\partial_{d}f( 0) )$. The vector
$\widetilde{W}_{u}$
is a $({d(d+1)}/{2})$-dimensional random vector with density function
\[
\bar{\psi}_{u}(w)\propto\bigl|{\det}\bigl(\overline{r}^{\ast}(w)-u\bar
{\Lambda}%
\bigr)\bigr|\exp\bigl(-\tfrac{1}{2}w^{\top}\mu_{2\cdot0}^{-1}w\bigr)\indic(\bar
{w}^{\ast}-u\bar{\Lambda
}\in\mathcal{N}),
\]
where $\bar{\Lambda}$ is the second spectral moment of $f$ restricted to
$\partial T$, and $\overline{r}{}^{\ast}(w)$ is the $(d-1)\times(d-1)$
symmetric
matrix whose upper triangular elements consist of the components of
$w$. In
the representation~(\ref{SlepianB}) the vectors~$\widetilde{W}_{u}$
and $Z$ are
independent. As in the proof of Lemma~\ref{LemBeta}, one can show that
a similar bound holds, albeit with with different constants. Thus, since
$\mu_{1}( t) =O( t) $, there exist
$c^{\prime\prime
}$ and $\delta^{\prime\prime}$ such that the third term in (\ref
{SlepianB}) can
be bounded by
\[
P\Bigl[\max_{\vert t\vert\leq\delta a/b}\vert\mu
_{1}(
t) \Lambda^{-1}( Z,0,\ldots,0) ^{T}\vert\geq
a/(4b)\Bigr]\leq c^{\prime\prime}\exp\biggl( -\frac{\delta^{\prime
\prime}}%
{\delta^{2}}\biggr) .
\]
Consequently, we can also find $c^{\ast}$ and $\delta^{\ast}$ such
that the
conclusion holds, and we are done.
%. For cases that the local maximum lies on lower dimensional
%boundaries, the proof is completely analogous.
\end{pf*}
\begin{pf*}{Proof of Lemma~\ref{LemGlobalMaxima}}
Recall that $t^*$ is the unique global maximum of~$f$ in~$T$.
Writing Palm probabilities as a ratio of expectations,
as explained in Section~\ref{extentpalmsec},
and using the fact that $t^*\in\mathcal L$, we immediately have
%
%e7.16 ###
%
\begin{eqnarray} \label{palmratio}
&& P\Bigl(\min_{\vert t-t^{\ast}\vert<\delta
ab^{-1}}f(
t) <b \big\|_{\mathcal P} f( t^{\ast})
>(b+a/b)\Bigr) \nonumber\\[-8pt]\\[-8pt]
&&\qquad \leq
\frac{ E(\#\{s\in\mathcal{L} \dvtx\min_{\vert
t-s\vert<\delta
ab^{-1}}f( t) <b , f( s) >b+a/b\})
}{E(\#\{s\in\mathcal L\dvtx f(s)>(b+a/b), s=t^* \}) }.\nonumber
\end{eqnarray}
Writing
\[
N_b=\#\{s\in\mathcal L\dvtx f(s)>(b+a/b)\},
\]
it is standard fare that,
for the random fields of the kind we are treating,
\[
E(N_b) = \bigl(1+o(1)\bigr)P(N_b=1)
\]
for large $b$; for example, Chapter 6 of~\cite{GRF} or Chapter 5 of
\cite{ADT10}.

Therefore, for $b$ large enough,
\[
\frac{E(\#\{s\in\mathcal L\dvtx f(s)>(b+a/b)\})}{E(\#\{s\in\mathcal
L\dvtx
f(s)>(b+a/b), s=t^* \})} < 2.
\]
Substituting this into~(\ref{palmratio}) yields, for any $s\in
\mathcal L$,
\begin{eqnarray*}
&&
P\Bigl(\min_{|t-t^{\ast}|<\delta ab^{-1}}f(t)<b \big\|_{\mathcal{P}}
f(t^{\ast
})>b+a/b\Bigr)
\\
&&\qquad
\leq2 P\Bigl(\min_{|t-s|<\delta ab^{-1}}f(t)<b \big\|_{\mathcal{P}}
f(s) > b+a/b\Bigr)
\\
&&\qquad \leq2c^{\ast}\exp( -\delta^{\ast}/\delta
^{2}),
\end{eqnarray*}
where the second inequality follows from~(\ref{Palm}), and we are done.
\end{pf*}

We complete the paper with the proofs of Lemmas~\ref{LemBeta} and~\ref{LemG}.
\begin{pf*}{Proof of Lemma~\ref{LemBeta}}
It suffices to prove the lemma for
the case $a=1$. Since
\[
f_{u}( 0) =u=u-W_{u}\beta^{\top}( 0)
+g(
0) ,
\]
and $W_{u}$ and $g $ are independent, $\beta(
0) =0$. Furthermore, since
$C^{\prime}( t) =O( t) $ and
$\mu_{2}^{\prime}( t) =O( t)$,
there exists a $c_{0}$ such that $|\beta(t)|\leq c_{0}|t|^{2}$. In addition,
$W_{u}$ has density function proportional to
\[
\psi_{u}(w)\propto\bigl|{\det}\bigl(r^{\ast}(w)-u\Lambda\bigr)\bigr|
\exp\bigl(-\tfrac
{1}{2}w^{\top}%
\mu_{2\cdot0}^{-1}w\bigr)\indic( w^{\ast}-u\Lambda\in\mathcal
{N}) .
\]
Note that $\det(r^{\ast}(w)-u\Lambda)$ is expressible as
a polynomial in $w$ and $u$,
and there exists some $\varepsilon_{0}$ and $c$ such that
\[
\biggl\vert\frac{\det( r^{\ast}(w)-u\Lambda) }{\det
(
-u\Lambda) }\biggr\vert\leq c,
\]
if $\vert w\vert\leq\varepsilon_{0}u$. Hence, there
exist $\varepsilon
_{2}$, $c_{2}>0$, such that
\[
\psi_{u}( w) \leq\widetilde{\psi}( w)
\dvtx=c_{2}\exp\bigl(
-\tfrac{1}{2}\varepsilon_{2}w^{\top}\mu_{2\cdot0}^{-1}w\bigr)
\]
for all $u\geq1$. The right-hand side here is proportional to a
multivariate Gaussian density. Thus,
\[
P( \vert W_{u}\vert>x) =\int_{\vert
w\vert
>x}\psi_{u}(w)\,dw\leq\int_{\vert w\vert>x}\widetilde
{\psi}_{u}(
w) \,dw=c_{3}P(|\widetilde{W}|>x),
\]
where $\widetilde{W}$ is a multivariate Gaussian random variable with density
function proportional to $\widetilde{\psi}$. Therefore, by choosing
$\varepsilon
_{1}$ and $c_{1}$ appropriately, we have
\[
P\biggl( \sup_{\vert t\vert\leq\delta/b}\vert
W_{u}%
\beta^{\top}\vert>\frac{1}{4b}\biggr) \leq P\biggl(
\vert
W_{u}\vert>\frac{b}{c_{0}^{2}\delta^{2}}\biggr) \leq
c_{1}\exp\biggl(
-\frac{\varepsilon_{1}b^{2}}{\delta^{4}}\biggr)
\]
for all $u\geq b$.
\end{pf*}
\begin{pf*}{Proof of Lemma~\ref{LemG}}
Once again, it suffices to prove
the lemma for
the case $a=1$. Since
\[
f_{b}(0)=b=b-W_{b}\beta^{\top}(0)+g(0),
\]
the covariance function $( \gamma( s,t) \dvtx
s,t\in T) $
of the centered field $g $ satisfies $\gamma(0,0)=0$.
It is also easy to check that
\[
\partial_{s}\gamma(s,t)=O(|s|+|t|),\qquad \partial_{t}\gamma(s,t)=O(|s|+|t|).
\]
Consequently, there exists a constant $c_{\gamma}\in(0,\infty)$ for which
\[
\gamma(s,t)\leq c_{\gamma}(|s|^{2}+|t|^{2}),\qquad \gamma(s,s)\leq
c_{\gamma
}|s|^{2}.
\]
We need to control the tail probability of $\sup_{|t|\leq\delta
/b}|g(t)|$. For
this it is useful to introduce the following scaling. Define
\[
g_{\delta}( t) = \frac
{b}{\delta} g\biggl( \frac{\delta t}{b}\biggr).
\]
Then $\sup_{\vert t\vert\leq\delta/b}g(
t) \geq\frac{1}{4b}$ if and only if $\sup_{\vert
t\vert
\leq1}g_{\delta}( t) \geq\frac{1}{4\delta}$. Let
\[
\sigma_{\delta}( s,t) =E( g_{\delta}(
s)
,g_{\delta}( t) ) .
\]
Then,
\[
\sup_{s\in\real}\sigma_{\delta}(s,s)\leq c_{\gamma}.%
\]
Because $\gamma(s,t)$ is at least twice differentiable, applying a Taylor
expansion we easily see that the canonical metric $d_g$ corresponding to
$g_{\delta}( s) $ (cf. Theorem~\ref{ThmA1}) can be
bounded as follows:
\begin{eqnarray*}
d_g^2(s,t) &=& E\bigl( g_{\delta}( s) -g_{\delta}(
t) \bigr)
^{2} \\
& = &\frac{b^{2}}{\delta^{2}}\biggl[ \gamma\biggl( \frac{\delta
s}%
{b},\frac{\delta s}{b}\biggr) +\gamma\biggl( \frac{\delta
t}{b},\frac{\delta
t}{b}\biggr) -2\gamma\biggl( \frac{\delta s}{b},\frac{\delta
t}{b}\biggr)
\biggr] \\
& \leq & c\vert s-t\vert^{2}
\end{eqnarray*}
for some constant $c\in(0,\infty$). Therefore, the entropy of
$g_{\delta
}$, evaluated at $\widetilde{\delta}$,
% namely $\mathcal{N}(\widetilde{\delta})$
%defined in Theorem~\ref{ThmA1} ,
is bounded by $K\widetilde{\delta}^{-d}$ for any
$\widetilde{\delta}>0$ and with an appropriate choice of $K>0$.
Therefore, for all
$\delta<\delta_{0}$,
\[
P\biggl(\sup_{|t|\leq\delta/b}|g(t)|\geq\frac{1}{4b}\biggr)=P\biggl(\sup_{
\vert
t\vert\leq1}g_{\delta}( t) \geq\frac
{1}{4\delta}\biggr)\leq
c_{d}\delta^{-d-\eta}\exp\biggl(-\frac{1}{16c_{\gamma}\delta^{2}}\biggr)
\]
for some constant $c_{d}$ and $\eta>0$. The last inequality is a direct
application of Theorem 4.1.1 of~\cite{AdlerTaylor07}. The conclusion
of the
lemma follows immediately by choosing $\widetilde{c}$ and $\widetilde
{\delta}$ appropriately.
\end{pf*}

%s8 ###
\section{Numerical examples}\label{SecNum}

In this section, we provide four examples which indicate how well the techniques
we have suggested actually work in practice.

The fist treats a random field for which the tail probability is in a
closed form. This is simply to confirm that the estimates yielded from
the algorithm are reasonable.
\begin{Example}
Let $f( t) =X\cos t+Y\sin t$ and $T=[ 0,3/4]
$ where $X$ and~$Y$ are i.i.d. standard Gaussian. We compute $P
(\sup_T f(t) >b) $. This probability is known in closed form
(cf.~\cite{AdlerTaylor07}) and is
given by
%
%e8.1 ###
%
\begin{equation} \label{CosExcur}
P\Bigl( \sup_{0\leq t \leq3/4} f(t) >b\Bigr) =1-\Phi(
b) +\frac{3}{%
8\pi}e^{-b^{2}/2}.
\end{equation}
Table~\ref{ExKl} shows the (remarkably accurate) simulation results.
\end{Example}

%
%t1 ###
%
\begin{table}
\caption{Simulation results for the cosine process. All results are
based on $10^{3}$ independent simulations. The ``True value'' is
computed using (\protect\ref{CosExcur})
The computation time for~each~estimate is less than one second. The
lattice size is $3b$}\label{ExKl}
\begin{tabular*}{\tablewidth}{@{\extracolsep{\fill}}lccc@{}}
\hline
$\bolds{b}$ & \textbf{True value} & \textbf{Est.} & \textbf{Std. er.} \\
\hline
%1.5 & 0.1055 & 0.1057 & 7.96E--04 \\
%2.0 & 0.0389 & 0.0392 & 4.02E--04 \\
%2.5 & 0.0115 & 0.0118 & 1.57E--04 \\
\hphantom{0}3 & 3.12E--03 & 3.13E--03 & 8.43E--05 \\
\hphantom{0}5 & 8.8E--07\hphantom{0} & 8.6E--07\hphantom{0} & 2.27E--08 \\
10 & 3.83E--23 & 3.81E--23 & 8.88E--25 \\
\hline
\end{tabular*}
\end{table}

The remaining examples treat more interesting random fields for
which~$T$ is a two-dimensional square.
\begin{Example}\label{ExHomo}
Consider the smooth homogenous random field on $T=[0,1]^2$ with mean
zero and covariance function
\[
C(t) = e^{-|t|^2}.
\]
Table~\ref{TabHomo} shows the simulation results of the excursion
%
%t2 ###
%
\begin{table}
\caption{Simulation results in Example \protect\ref{ExHomo}}\label{TabHomo}
\begin{tabular*}{\tablewidth}{@{\extracolsep{\fill}}lcccccc@{}}
\hline
& \multicolumn{2}{c}{$\bolds{P(\sup_{T}f(t)>b)}$} & \multicolumn
{2}{c}{$\bolds{E(
\sup_{T}f(t)-b |{\sup_{T}f}(t)>b)}$} & & \\[-4pt]
& \multicolumn{2}{c}{\hrulefill} & \multicolumn
{2}{c}{\rule{124pt}{0.5pt}} & &\\
$\bolds{b}$ & \textbf{Est.} & \textbf{St. d.}
& \multicolumn{1}{c}{\hspace*{23pt}\textbf{Est.}} &
\multicolumn{1}{c}{\textbf{St. d.}} & \textbf{Lattice size} & \textbf{CPU time} \\
\hline
3 & 1.1E--02 & 3.8E--04 & \hspace*{23pt}0.30 & 1.5E--02 & 10 by 10 & \hphantom{0}6 sec \\
4 & 3.3E--04 & 1.2E--05 & \hspace*{23pt}0.25 & 1.3E--02 & 15 by 15 & 53 sec \\
5 & 4.3E--06 & 1.6E--07 &  \hspace*{23pt}0.19 & 1.0E--02 & 15 by 15 & 45 sec \\
%8 & 1.9E--14 & 7.2E--16 & & & 2.4E--15 & & 5.0E--17 & & \\
\hline
\end{tabular*}  \vspace*{-12pt}
\end{table}
probabilities\break $P(\sup_T f(t) >b)$ and expected overshoots $E
(\sup_T f(t) -b | {\sup_T f}(t) >b)$. The results are based
on 1\mbox{,}000 independent simulations by setting the tuning parameter $a=1$.
The size of discretization and CPU time are also reported.
\end{Example}

\begin{Example}\label{ExVarMean}
Consider the continuous, but nondifferentiable, and nonhomogenous
random field on $T=[0,1]^2$ with
\[
\mu(t) = 0.1 t_1 + 0.1 t_2 C(s,t) = e^{-|t-s|^2}.
\]
Table~\ref{TabVarMean} shows the simulation results of excursion
probabilities $P(\sup_T f(t) >b)$ and expected overshoots
%
%t3 ###
%
\begin{table}
\caption{Simulation results for Example \protect\ref{ExVarMean}}\label{TabVarMean}
\begin{tabular*}{\tablewidth}{@{\extracolsep{\fill}}lcccccc@{}}
\hline
& \multicolumn{2}{c}{$\bolds{P(\sup_{T}f(t)>b)}$} & \multicolumn
{2}{c}{$\bolds{E(
\sup_{T}f(t)-b |{\sup_{T}f}(t)>b)}$} & & \\[-4pt]
& \multicolumn{2}{c}{\hrulefill} & \multicolumn
{2}{c}{\rule{124pt}{0.5pt}} & &\\
$\bolds{b}$ & \textbf{Est.} & \textbf{St. d.}
& \multicolumn{1}{c}{\hspace*{23pt}\textbf{Est.}} &
\multicolumn{1}{c}{\textbf{St. d.}} & \textbf{Lattice size} & \textbf{CPU time} \\
\hline
3 & 1.4E--02 & 5.0E--04 & \hspace*{23pt}0.32 & 1.6E--02 & 10 by 10 & \hphantom{0}6 sec \\
4 & 5.3E--04 & 1.9E--05 & \hspace*{23pt}0.25 & 1.3E--02 & 15 by 15 & 40 sec \\
5 & 7.2E--06 & 2.6E--07 & \hspace*{23pt}0.20 & 9.8E--03 & 15 by 15 & 56 sec \\
%8 & 4.4E--14 & 1.4E--15 & & & 5.6E--15 & & 1.2E--16 & & \\
\hline
\end{tabular*}
\end{table}
$E(\sup_T f(t) -b | {\sup_T f}(t) >b)$. The simulation
setting is the same as that in Example~\ref{ExHomo}.
\end{Example}
\begin{Example}\label{ExVarMeanHol}
Consider the smooth random field living on $T=[0,1]^2$ with
\[
\mu(t) = 0.1 t_1 + 0.1 t_2 C(t) = e^{-|t|/4}.
\]
Table~\ref{TabVarMeanHol} shows simulation results for the excursion
probabilities $P(\sup_T f(t) >b)$ and the expected overshoots
%
%t4 ###
%
\begin{table}[b]
\caption{Simulation results for Example \protect\ref{ExVarMeanHol}}\label{TabVarMeanHol}
\begin{tabular*}{\tablewidth}{@{\extracolsep{\fill}}lcccccc@{}}
\hline
& \multicolumn{2}{c}{$\bolds{P(\sup_{T}f(t)>b)}$} & \multicolumn
{2}{c}{$\bolds{E(
\sup_{T}f(t)-b |{\sup_{T}f}(t)>b)}$} & & \\[-4pt]
& \multicolumn{2}{c}{\hrulefill} & \multicolumn
{2}{c}{\rule{124pt}{0.5pt}} & &\\
$\bolds{b}$ & \textbf{Est.} & \textbf{St. d.}
& \multicolumn{1}{c}{\hspace*{23pt}\textbf{Est.}} &
\multicolumn{1}{c}{\textbf{St. d.}} & \textbf{Lattice size} & \textbf{CPU time} \\
\hline
3 & 1.5E--02 & 5.8E--04 & \hspace*{23pt}0.33 & 1.5E--02 & 15 by 15 & \hphantom{0}58 sec \\
4 & 6.4E--04 & 3.1E--05 & \hspace*{23pt}0.25 & 1.4E--02 & 15 by 15 & \hphantom{0}44 sec \\
5 & 1.3E--05 & 6.9E--07 & \hspace*{23pt}0.21 & 1.3E--02 & 25 by 25 & 600 sec\\
%7 & 1.2E--09 & 6.7E--11 & & & 1.9E--10 & & 5.9E--12 & &
\hline
\end{tabular*}
\end{table}
$E(\sup_T f(t) -b | {\sup_T f}(t) >b)$. The simulation setting is
the same as that in Example~\ref{ExHomo}.
\end{Example}

Although we have given rigorous results regarding descretization
parameters, in
practice we choose the lattice size sufficiently large so that the bias
was inconsequential in comparison to the estimated standard deviation.\vadjust{\goodbreak}
We achieved this by increasing the lattice size until the change of the
estimate was small enough relative to the estimated standard
deviation.

Note that, for all the examples, the relative error does not increase
as the
level increases and the exceedance
probability tends to zero as long as the lattice size also increases.
This is in line with the theoretical results of the paper.\looseness=-1

Another empirical finding is that the computational burden increases
substantially with lattice size, although the algorithm has been proven
to be of polynomial complexity. This complexity is mainly from the
Cholesky decomposition of large covariance matrices.
While this is a problem common to all discrete simulation algorithms
for random fields, we nevertheless plan to look at this efficiency
issue in future
work.

\section*{Acknowledgment}
We are grateful to a referee and an Associate Editor for helpful
comments and suggestions.

%suskaldyti doi

% imsref loaded by lrinkeviciute, 2011-08-09 12:54:08
%
% imsref loaded by lrinkeviciute, 2011-08-10 09:54:50
% imsref loaded by lrinkeviciute, 2011-08-10 09:55:38
% imsref loaded by lrinkeviciute, 2011-08-10 09:58:08
% imsref loaded by lrinkeviciute, 2011-08-10 10:54:13

%
\printaddresses

\end{document}